\documentclass[11pt,british,reqno]{amsart}

\usepackage[T1]{fontenc}

\usepackage{babel,amssymb,booktabs,eucal,paralist,csquotes,stmaryrd}
\usepackage[margin=30pt,font=small]{caption}

\usepackage{ifpdf}
\ifpdf
  \usepackage{microtype}
  \usepackage[bookmarks,bookmarksnumbered=true]{hyperref}
\fi

\theoremstyle{plain}
\newtheorem{theorem}{Theorem}[section]
\newtheorem{proposition}[theorem]{Proposition}
\newtheorem{lemma}[theorem]{Lemma}
\newtheorem{corollary}[theorem]{Corollary}

\theoremstyle{definition}
\newtheorem{definition}[theorem]{Definition}

\theoremstyle{remark}

\newtheorem{remark}[theorem]{Remark}
\newtheorem*{acknowledgements}{Acknowledgements}

\numberwithin{equation}{section}
\numberwithin{table}{section}

\setdefaultenum{\upshape (i)}{\upshape (a)}{}{}

\newcommand{\Lie}[1]{\operatorname{\textsl{#1}}}
\newcommand{\lie}[1]{\operatorname{\mathfrak{#1}}}
\newcommand{\GL}{\Lie{GL}}
\newcommand{\SO}{\Lie{SO}}
\newcommand{\so}{\lie{so}}
\newcommand{\SP}{\Lie{Sp}}
\newcommand{\sP}{\lie{sp}}
\newcommand{\SU}{\Lie{SU}}
\newcommand{\Un}{\Lie{U}}
\newcommand{\HQ}{\Lie H_Q}
\newcommand{\Gk}{\Lie G(k)}

\newcommand{\LC}{\nabla^{\textup{LC}}}
\newcommand{\Nt}{\tilde\nabla}
\newcommand{\NaqH}{\nabla^{\textup{aqH}}}
\newcommand{\Nhkt}{\nabla^{\textup{HKT}}}
\newcommand{\Nqkt}{\nabla^{\textup{QKT}}}
\newcommand{\Nq}{\nabla^{\textup{q}}}
\newcommand{\IaH}{\xi^{\textup{aH}}}

\DeclareMathOperator{\End}{End}
\DeclareMathOperator{\Alt}{\textit{Alt}}

\newcommand{\hook}{\mathbin{\lrcorner}}
\newcommand{\inp}[2]{\left< #1, #2\right>}
\newcommand{\norm}[1]{\left\lVert #1\right\rVert}

\newcommand{\sym}[1]{{S^{#1}\mskip-2mu}}

\newcommand{\Wc}[1]{\mathcal W_{#1}}
\DeclareMathOperator{\Wm}{W}

\DeclareFontFamily{U}{bigeuf}{}
\DeclareFontShape{U}{bigeuf}{m}{n}{<-6>s*[1.5]eufm5%
<6-8>s*[1.5]eufm7%
<8->s*[1.5]eufm10}{}
\DeclareSymbolFont{bigeufletters}{U}{bigeuf}{m}{n}
\DeclareMathSymbol{\sumcic}{\mathop}{bigeufletters}{`S}

\newcommand{\bdash}{-\hspace{0pt}} 
\newcommand{\eqbreak}{\\&\qquad}
\newcommand{\smalleqbreak}{\\&\quad}
\newcommand{\detc}{\quad\text{etc.}}

\newcommand{\itref}[1]{\textup{(\ref{#1})}}
\makeatletter
\renewcommand{\p@enumii}{\theenumi.}
\makeatother

\begin{document}

\title[Intrinsic torsion of almost quaternion-Hermitian manifolds]%
{The intrinsic torsion of almost quaternion-Hermitian manifolds}

\author{Francisco Mart\'\i n Cabrera}
\address[F.~Mart\'\i n Cabrera]{Department of Fundamental Mathematics\\
University of La Laguna\\ 38200 La Laguna, Tenerife, Spain}
\email{fmartin@ull.es}

\author{Andrew Swann} \address[A.F.~Swann]{Department of Mathematics and
Computer Science\\ University of Southern Denmark\\ Campusvej 55, DK-5230
Odense M, Denmark} \email{swann@imada.sdu.dk}

\begin{abstract}
  We study the intrinsic torsion of almost quaternion\bdash Hermitian manifolds
  via the exterior algebra.  In particular, we show how it is determined by
  particular three-forms formed from simple combinations of the exterior
  derivatives of the local K\"ahler forms.  This gives a practical method
  to compute the intrinsic torsion and is applied in a number of examples.
  In addition we find simple characterisations of HKT and QKT geometries
  entirely in the exterior algebra and compute how the intrinsic torsion
  changes under a twist construction.
\end{abstract}



\maketitle

\begin{center}
  \begin{minipage}{0.7\linewidth}
    \begin{small}
      \tableofcontents
    \end{small}
  \end{minipage}
\end{center}
\section{Introduction}
\label{sec:introduction}

An almost quaternion-Hermitian manifold~$M$ is a Riemannian $4n$-manifold
which admits an $\SP(n)\SP(1)$-structure, i.e., a reduction of its frame
bundle to the subgroup $\SP(n)\SP(1)$ of~$\SO(4n)$.  This is equivalent to
the presence of a Riemannian metric $g=\inp\cdot\cdot$ and a rank-three
subbundle $\mathcal G$ of the endomorphism bundle~$\End TM$, locally
generated by three almost complex structures $I$, $J$, $K$ satisfying the
identities of the imaginary unit quaternions.  Almost quaternion-Hermitian
manifold are of special interest because $\SP(n)\SP(1)$ is included in
Berger's list~\cite{Berger:hol} of possible holonomy groups of locally
irreducible Riemannian manifolds that are not locally symmetric.  Also in
the field of theoretical physics, the study of supersymmetric sigma models
and their couplings to supergravity is very related with the study of
complex and quaternionic structures defined on Riemannian
manifolds~\cite{Gates-HR:twisted,Howe-P:further}.

Since $\SP(n)\SP(1)$ is a closed and connected subgroup of~$\SO(4n)$, there
exists a unique metric $\SP(n)\SP(1)$-connection $\NaqH = \LC + \xi$, where
$\LC$ is the Levi-Civita connection and $\xi$ is a tensor, called the
\emph{intrinsic $\SP(n)\SP(1)$-torsion}, in $T^*M \otimes
(\sP(n)+\sP(1))^\perp$.  Here $(\sP(n)+\sP(1))^\perp$ denotes the
orthogonal complement in~$\so(4n)$ of the Lie algebra~$\sP(n)+\sP(1)$.

Under the action of $\SP(n)\SP(1)$, the space $T^*M \otimes
(\sP(n)+\sP(1))^{\perp}$ of possible intrinsic torsion tensors $\xi$
decomposes into irreducible $\SP(n)\SP(1)$-modules, giving rise to a
natural classification of almost quaternion\bdash Hermitian manifolds.
In~\cite{Swann:symplectiques} it was shown that, in general dimensions,
$\xi$ has six components and $2^6=64$ classes of such manifolds potentially
arise.  An almost quaternion-Hermitian manifold is said to be
quaternion-K\"ahler, if the intrinsic torsion $\xi$ vanishes.  In this
case, the reduced holonomy group is a subgroup of $\SP(n)\SP(1)$ and the
manifold is Einstein.  On the other hand, if the three almost complex
structures are globally defined, then $M$ is said to be endowed with an
almost hyperHermitian structure (an $\SP(n)$-structure).  When the three
K\"ahler forms $\omega_I$, $\omega_J$, $\omega_K$ of the $\SP(n)$-structure
are covariant constant, the manifold is called hyperK\"ahler.
HyperK\"ahler manifolds have reduced holonomy group contained in $\SP(n)$
and their Ricci curvature vanishes.

By identifying the intrinsic $\SP(n)\SP(1)$-torsion~$\xi$ with the
Levi-Civita covariant derivative of a certain four-form~$\Omega$, defined
below in equation~\eqref{eq:four-form}, one obtains an analogue of the
method of Gray \& Hervella \cite{Gray-H:16} for finding conditions for
classes of almost quaternion-Hermitian manifolds.  Detailed conditions
describing classes in this way were given in~\cite{Cabrera:aqh}.

In the present paper, we will take another approach.  In fact, we will show
how the intrinsic torsion~$\xi$ can be determined by means of the exterior
derivatives $d \omega_I$, $d \omega_J$ and $d \omega_K$ of the local
K\"ahler forms corresponding to the almost complex structures $I$,
$J$,~$K$.  In the process, there will arise additional, detailed
information about the components of~$\xi$ which will be very useful in
working on examples of almost quaternion-Hermitian manifolds.  For all of
this, we give expressions for the covariant derivatives $\LC \omega_A$ in
terms of $d \omega_I$, $d \omega_J$ and $d \omega_K$, see
Proposition~\ref{prop:naddd}.  Such expressions contribute to a better
understanding of Hitchin's result \cite{Hitchin:Riemann-surface} saying
that if $\omega_I$, $\omega_J$ and $\omega_K$ are closed, then they are
covariant constant.  Indeed in Proposition~\ref{prop:naddd}, we show how
the Nijenhuis tensor~$N_I$ in general is determined by the difference $J d
\omega_J - K d \omega_K$.  Let us briefly explain one application of this
result, cf.~\S\ref{sec:HKT}.

It is known that the geometry of the target space of $(4,0)$ supersymmetric
even-dimensional sigma models without Wess-Zumino term (torsion) is a
hyperK\"ahler manifold.  In presence of torsion, the geometry of the target
space is a hyperK\"ahler manifold with torsion, usually called an
HKT-manifold \cite{Howe-P:twistor-kaehler}.

Grancharov \& Poon \cite{Grantcharov-P:HKT} showed that an almost
hyperHermitian manifold $(M,I,J,K,\inp\cdot\cdot)$ is HKT if and only if:
\begin{enumerate}
\item\label{item:integrable} the three almost complex structures $I$, $J$
  and $K$ are integrable, and
\item\label{item:Id} $Id \omega_I= Jd \omega_J = K d \omega_K$.
\end{enumerate}
A direct consequence of our expression for~$N_A$ is that
condition~\itref{item:Id} is sufficient to characterise HKT geometry, and
in particular \itref{item:Id} implies the integrability
condition~\itref{item:integrable}.  Similarly, we also show how Grancharov
\& Poon's holomorphic characterisation for HKT-manifolds may be simplified,
see~\S\ref{sec:HKT}.

It is known that an almost quaternion-Hermitian manifold the three
covariant derivatives $\LC \omega_I$, $\LC \omega_J$ and $\LC \omega_K$ are
not independent, but rather any two determine the third
\cite{Fernandez-Moreiras:symmetry,Cabrera:aqh}, see
equation~\eqref{eq:sym-nabla}.  The corresponding statement for the
exterior derivatives $d \omega_I$, $d \omega_J$, $d \omega_K$ is not true.
However, we find that there are still relations expressed by symmetries of
the three-forms
\begin{equation*}
  \beta_I = J d \omega_J + K d \omega_K, \detc
\end{equation*}
see~\eqref{eq:naddd4r}.  These symmetries are equivalent to requiring
$\beta_A$ to be of type $\{2,1\}=(2,1)+(1,2)$ with respect to the almost
complex structure~$A$.  Algebraically the $\beta_I$, $\beta_J$ and
$\beta_K$ are independent three-forms of these types and we find that the
space of possible triples of covariant derivatives $(\nabla \omega_I,
\nabla \omega_J, \nabla \omega_K)$ is isomorphic to the space of possible
triples of three-forms $(\beta_I, \beta_J, \beta_K)$.

The relevance of the three-forms $\beta_I$, $\beta_J$, $\beta_K$ is clearly
seen in Proposition~\ref{prop:torsion-beta}, where we demonstrate how they
may be used to compute the components of the intrinsic
$\SP(n)\SP(1)$-torsion~$\xi$.  This gives a practical way to compute~$\xi$
via the exterior algebra and will be used in the study of concrete examples
in~\S\ref{sec:examples}.

In~\S\ref{sec:QKT}, we focus attention on quaternion-K\"ahler manifolds
with torsion, also known as QKT-manifolds.  Motivation for studying these
structures can be also found in the field of theory of supersymmetric sigma
models, see~\cite{Howe-OP:QKT}.  Our results lead to a new
characterisation~\eqref{eq:omega-diff} of QKT-manifolds that is simpler
than that provided by Ivanov \cite[Theorem~2.2]{Ivanov:QKT}.  We also
obtain new expressions for the torsion three-form and torsion one-form and
study of the integrability properties of the almost complex structures.

In~\S\ref{sec:twist}, we consider the intrinsic torsion of
quaternion-Hermitian manifolds obtained by the twist interpretation of
T-duality given in~\cite{Swann:T}.  Using the exterior algebra is
particularly advantageous here.  We see that in many cases the QKT
condition is preserved.

Finally, in~\S\ref{sec:examples}, we give an number of examples of types of
almost quaternion\bdash Hermitian manifolds.  We particularly mention one
of the quaternionic structures considered on the manifold $S^3 \times T^9$
which is a non-QKT-manifold admitting an $\SP(n) \SP(1)$-connection with
skew-symmetric torsion $(0,3)$-tensor.  In the various examples, we have
also determined the types of almost Hermitian structure.  Because it is an
advantage to handle Lie brackets instead of directly using~$\LC$, we
determine such types by means of the exterior derivative~$d \omega_I$ and
the Nijenhuis tensor~$N_I$.  Therefore, in \S\ref{sec:Gray-Hervella}, we
include Table~\ref{tab:ahtypes} showing conditions in terms of $d \omega_I$
and $N_I$ to characterise the Gray-Hervella classes.

\begin{acknowledgements}
  This work is supported by a grant from the MEC (Spain), project no.\
  MTM2004-2644.  We thank the organisers of the Workshop on
  \enquote{Special Geometries in Mathematical Physics}, K\"uhlungsborn,
  2006, for the chance to present some of this material.
\end{acknowledgements}

\section{Definitions and notation}
\label{sec:notation}

Let $G$ be a subgroup of the linear group $\GL(m, \mathbb R)$.  A
manifold~$M$ is said to be equipped with a $G$-structure, if there is a
principal $G$-subbundle~$P$ of the principal frame bundle.  In such a case,
there always exist connections, called \emph{$G$-connections}, defined on
the subbundle~$P$.  Moreover, if $(M^m ,\inp\cdot\cdot)$ is an orientable
$m$-dimensional Riemannian manifold and $G$ a closed and connected subgroup
of~$\SO(m)$, then there exists a unique metric $G$-connection $\nabla^G$
such that $\xi = \nabla^G - \LC$ takes its values in $\lie g^\perp$, where
$\lie g^\perp$ denotes the orthogonal complement in~$\so(m)$ of the Lie
algebra~$\lie g$ of~$G$ \cite{Salamon:holonomy,Cleyton-S:intrinsic}.  The
tensor $\xi$ is said to be the \emph{intrinsic $G$-torsion} and $\nabla^G$
is called the \emph{minimal $G$-connection}.

A $4n$-dimensional manifold $M$ is said to be \emph{almost
quaternion-Hermitian}, if $M$ is equipped with an $\SP(n)\SP(1)$-structure.
This is equivalent to the presence of a Riemannian metric $\inp\cdot\cdot$
and a rank-three subbundle $\mathcal G$ of the endomorphism bundle~$\End
TM$, such that locally $\mathcal G$ has an \emph{adapted basis} $I, J, K$
satisfying $I^2 = J^2= -1$ and $K= IJ = - JI$, and $\inp{AX}{AY} = \inp XY
$, for all $X,Y \in T_x M$ and $A =I,J,K$.  An almost quaternion-Hermitian
manifold with a global adapted basis is called an \emph{almost
hyperHermitian} manifold.  In such a case the structure group reduces
to~$\SP(n)$.  We note that if $I,J,K$ is an adapted basis then so are
$J,K,I$ and $K,I,J$; thus formul\ae\ derived for an arbitrary adapted
$I,J,K$ will also apply to cyclic permutations of these almost complex
structures.

There are three local K\"ahler-forms $\omega_A (X,Y) = \inp X{AY}$,
$A=I,J,K$.  From these one may define a global, non-degenerate
four-form~$\Omega$, the \emph{fundamental form}, via the local formula
\begin{equation}
  \label{eq:four-form}
  \Omega = \sum_{A=I,J,K} \omega_A \wedge \omega_A.
\end{equation}
We will write
\begin{equation*}
  \Lambda_I\colon \Lambda^pT^*M \to \Lambda^{p-2}T^*M
\end{equation*}
for the adjoint of $\cdot\mapsto\cdot\wedge\omega_I$ with respect to the
metrics
\begin{equation*}
  \inp ab = \tfrac1{p!} a(e_{i_1}, \dots, e_{i_p})
  b(e_{i_1},\dots, e_{i_p}).
\end{equation*}
In particular, for a three-form~$\beta$ we have $\Lambda_I\beta =
\inp{\cdot\hook\beta}{\omega_I}$, and for a one-form $\nu$, we have
\begin{equation*}
  \Lambda_I(\nu\wedge\omega_A)=
  -\tfrac12\sum_{i=1}^{4n} \bigl( \nu(e_i)\omega_A(Ie_i,\cdot) +
  \omega_A(e_i,Ie_i)\nu  + \nu(Ie_i)\omega_A(\cdot,e_i) \bigr). 
\end{equation*}

In the next section, we will explicitly describe the intrinsic torsion of
almost quaternion\bdash Hermitian manifolds.  For such a purpose, we need
some basic tools related with almost quaternion-Hermitian manifolds in a
context of representation theory.  We will follow the $E$-$H$-formalism
used in \cite{Salamon:Invent,Swann:symplectiques} and we refer to
\cite{Broecker-tom-Dieck:Lie} for general information on representation
theory.  Thus, $E$~is the fundamental representation of $\SP(n)$ on
$\mathbb C^{2n} \cong \mathbb H^n$ via left multiplication by quaternionic
matrices, considered in $\GL(2n, \mathbb C)$, and $H$ is the representation
of $\SP(1)$ on $\mathbb C^2 \cong \mathbb H$ given by $q . \zeta = \zeta
\bar q$, for $q \in \SP(1)$ and $\zeta \in \mathbb H$.  An
$\SP(n)\SP(1)$-structure on a manifold $M$ gives rise to local bundles $E$
and $H$ associated to these representation and identifies $T M
\otimes_{\mathbb R} \mathbb C \cong E \otimes_{\mathbb C} H$.

On $E$, there is an $\SP(n)$-invariant complex symplectic form $\omega_E$
and a Hermitian inner product given by ${\inp xy}_{\mathbb C} = \omega_E(x,
\tilde y)$, where $y\mapsto\tilde y = j y $ is a quaternionic structure map
on $E=\mathbb C^{2n}$ considered as left complex vector space.  The mapping
$x \mapsto x^\omega = \omega_E ( \cdot, x)$ gives us an identification
of~$E$ with its dual~$E^*$.  If $\{ u_1, \dots , u_n, \tilde u_1, \dots ,
\tilde u_n \}$ is a complex orthonormal basis for~$E$, then $ \omega_E =
u^\omega_i \wedge \tilde u^\omega_i = u^\omega_i \tilde u^\omega_i - \tilde
u^\omega_i u^\omega_i$, where we have used the summation convention and
omitted tensor product signs.  These conventions will be used throughout
the paper.

The $\SP(1)$-module $H$ will be also considered as a left complex vector
space.  Regarding $H$ as a $4$-dimensional real space with the Euclidean
metric $\inp\cdot\cdot$ such that $\{1,i,j,k\}$ is an orthonormal basis.
The complex symplectic form is given by $\omega_H = 1^\flat \wedge j^\flat
+ k^\flat \wedge i^\flat + i (1^\flat \wedge k^\flat + i^\flat \wedge
j^\flat)$, where $h^\flat$ is given by $q \mapsto \inp hq$.  We also have
the identification, $h \mapsto h^\omega = \omega_H (\cdot , h )$, of $H$ with
its dual $H^*$ as complex space.  On $H$, we have a quaternionic structure
map given by $q = z_1 + z_2 j \mapsto \tilde q = jq = - \bar z_2 + \bar z_1
j$, where $z_1,z_2 \in \mathbb C$ and $\bar z_1, \bar z_2 $ are their
conjugates.  If $h \in H$ is such that $\inp hh =1$, then $\{ h , \tilde h
\}$ is a basis of the complex vector space $H$ and $\omega_H = h^\omega
\wedge \tilde h^\omega$.

The irreducible representations of $\SP(1)$ are the symmetric powers $\sym
k H \cong \mathbb C^{k+1}$.  An irreducible representation of $\SP(n)$ is
determined by its dominant weight $(\lambda_1,\dots,$ $ \lambda_n)$, where
$\lambda_i$ are integers with $\lambda_1 \geqslant \lambda_2 \geqslant
\dots \geqslant \lambda_n \geqslant 0$.  This representation will be
denoted by $V^{ (\lambda_1, \dots , \lambda_r)}$, where $r$ is the largest
integer such that $\lambda_r > 0$.  We will only need to use some of these
representations and use more familiar notation for these: $\sym k E =
V^{(k)}$, the $k$th symmetric power of~$E$; $\Lambda_0^r E = V^{(1, \dots ,
1)}$, where there are $r$ ones in exponent and $\Lambda_0^r E$ is the
$\SP(n)$-invariant complement to $\omega_E \Lambda^{r-2} E$ in $\Lambda^r
E$; also $K=V^{(21)}$, which arises in the decomposition $ E \otimes
\Lambda_0^2 E \cong \Lambda_0^3 E + K + E$, where $+$ denotes direct sum.

Most of the time in this paper, if $V$ is a complex $G$-module equipped
with a real structure, $V$ will also denote the real $G$-module which is
$(+1)$-eigenspace of the structure map.  The context should tell us which
space we are referring to.  However, when a risk of confusion arise, we
will denote the second mentioned space by $[V]$.  Likewise, the following
conventions will be used in this paper.  If $\psi$ is a $(0,s)$-tensor, for
$A=I,J,K$, we write
\begin{gather*}
    A_{(i)}\psi(X_1, \dots, X_i, \dots , X_s)
    = - \psi(X_1, \dots , AX_i, \dots , X_s),\\
    A_{(ij\dots k)} = A_{(i)}A_{(j)}\dots A_{(k)}, \quad\text{and}\quad \\
    A \psi(X_1,\dots,X_s) = (-1)^s \psi(AX_1,\dots,AX_s).
\end{gather*}

\section{The intrinsic torsion via differential forms}
\label{sec:intrinsic}

The intrinsic $\SP(n)\SP(1)$-torsion~$\xi$, $n>1$, is in $ T^*M \otimes
\left( \sP(n) + \sP(1) \right)^\perp \cong EH \otimes \Lambda_0^2 E \sym2 H
\subset T^*M \otimes \Lambda^2 T^*M$.  The space $EH \otimes \Lambda_0^2
E \sym2 H$ consists of tensors $\zeta$ such that
\begin{enumerate}
\item $(1 + I_{(23)} + J_{(23)} + K_{(23)})\zeta = 0$;
\item $\Lambda_A(\zeta_X) = 0$, for $A = I,J,K$, $X\in TM$,
\end{enumerate}
where $I,J,K$ is an adapted basis of~$\mathcal G$.  We recall that
$\Lambda^2 T^*M = \sym2 E + \sym2 H + \Lambda_0^2 E \sym2 H $, where $\sym2
E \cong \sP(n)$ and $\sym2 H \cong \sP(1)$ are the Lie algebras of $\SP(n)$
and $\SP(1)$, respectively.

A connection $\Nt$ is an $\SP(n)\SP(1)$-connection if $\Nt \Omega =0$.
This is the same as saying that $\Nt$ is metric, $\Nt g=0$, and
\emph{quaternionic}, meaning that for any local adapted basis $I$, $J$, $K$
of $\mathcal G$ we have
\begin{equation}
  \label{eq:q-nabla}
  (\Nt_X I) Y = \gamma_K(X) J Y - \gamma_J (X) KY, \detc,
\end{equation}
where $\gamma_I$, $\gamma_J$ and $\gamma_K$ are locally defined one-forms.
Here and throughout the rest of this paper, `etc.'\ means the equations
obtained by cyclically permuting $I,J,K$.

\begin{proposition}
  \label{prop:torsion}
  The minimal $\SP(n)\SP(1)$-connection is given by
  \begin{equation*}
    \NaqH = \LC + \xi,
  \end{equation*}
  where $\LC$ is the Levi-Civita connection and the intrinsic
  $\SP(n)\SP(1)$-torsion $\xi$ is given by
  \begin{equation*}
    \xi_X Y = - \tfrac14 \sum_{A=I,J,K} A (\LC_X A) Y + \tfrac12
    \sum_{A=I,J,K} \lambda_A(X) A Y,  
  \end{equation*}
  for all vectors $X,Y$.  Here the one-forms $\lambda_I$, $\lambda_J$ and
  $\lambda_K$ are defined by
  \begin{equation}
    \label{eq:lambda-I}
    \lambda_I(X) = \tfrac1{2n} \inp{\LC_X\omega_J}{\omega_K}, \detc 
  \end{equation}
\end{proposition}

\noindent
Note that if $n=1$, then $\NaqH = \LC$ and $\xi=0$.

\begin{proof}
  It is not hard to check $\NaqH g =0$, so $\NaqH$ is metric.  Now,
  computing $ (\xi_X I) Y = \xi_X IY - I \xi_X Y$, it is straightforward to
  obtain
  \begin{equation}
    \label{eq:aqH-torsion}
    (\LC_X I) Y = \lambda_K(X) JY - \lambda_J (X) KY - \xi_X IY + I \xi_X
    Y, \quad \text{etc}.
  \end{equation}
  Hence $(\NaqH_X I) Y = (\LC_X I) Y + (\xi_X I) Y= \lambda_K(X) JY -
  \lambda_J(X) KY$.  Therefore, $\NaqH$ is an $\SP(n)\SP(1)$-connection.
  
  Furthermore, the tensor $\xi$ satisfies
  \begin{equation*}
    \sum_{A=I,J,K} A \xi_X A Y = \xi_X Y
    \quad\text{and}\quad
    \sum_{i=1}^{4n} \inp{\xi_X e_i}{A e_i} =0,
  \end{equation*}
  for $A=I,J,K$.  Since these conditions imply $\xi \in T^*M \otimes
  (\sP(n) + \sP(1))^\perp = T^*M \otimes \Lambda_0^2 E \sym2H$, then
  $\NaqH = \LC + \xi$ is the minimal $\SP(n)\SP(1)$\bdash connection.
\end{proof}

The next result describes the decomposition of the space of possible
intrinsic torsion tensors $ T^*M \otimes \Lambda_0^2 E \sym2H$ into
irreducible $\SP(n)\SP(1)$-modules.

\begin{theorem}[Swann~\cite{Swann:symplectiques}]
  The intrinsic torsion $\xi$ of an almost quaternion-Hermitian manifold
  $M$ of dimension at least $8$, has the property
  \begin{equation}
    \label{eq:swann-sum}
    \begin{split}
      \xi \in T^*M \otimes \Lambda_0^2E \sym2 H &= \Lambda_0^3 E \sym3 H +
      K \sym3 H + E \sym3 H \eqbreak + \Lambda_0^3 E H + K H + E H.\qed
    \end{split}
  \end{equation}
\end{theorem}

If the dimension of $M$ is at least~$12$, all the modules of the sum are
non-zero.  For an eight-dimensional manifold $M$, we have $\Lambda_0^3 E =
\{ 0 \}$.  Therefore, for $\dim M \geqslant 12 $, we have $2^6 =64$ classes of
almost quaternion\bdash Hermitian manifolds, whereas there are $2^4=16$
classes when $\dim M =8$.  The map $\xi\mapsto \LC\Omega=-\xi\Omega$ is an
isomorphism, and in \cite{Cabrera:aqh} this was exploited to give explicit
conditions characterising these classes in terms of conditions
on~$\LC\Omega$.  However, from such conditions, it is not hard to derive
descriptions for the corresponding $\SP(n)\SP(1)$-components of~$\xi$ as we
will now demonstrate.

Firstly, the space of three-forms $\Lambda^3 T^*M$ decomposes under the
action of $\SP(n) \SP(1)$ as
\begin{equation*}
  \Lambda^3 T^*M = \Lambda^3_0 E \sym3 H +E \sym3 H + K H + EH.
\end{equation*}
Consider the operator
\begin{equation}
  \label{eq:L}
  \mathcal L = \mathcal L_I + \mathcal L_J + \mathcal L_K
\end{equation}
on $\Lambda^3T^*M$, where
\begin{equation*}
  \mathcal L_A = A_{(12)} + A_{(13)} + A_{(23)}.
\end{equation*}
The operator~$\mathcal L$ has eigenvalues $+3$ and $-3$ with corresponding
eigenspaces $(K+E)H$ and $(\Lambda^3_0E+E)\sym3H$.  For $\psi \in \Lambda^3
T^*M$, we have $\psi=\psi_H+\psi_{\sym3H}$ with
\begin{gather}
  \label{eq:psi-H}
  \psi_H =\tfrac16 (3\psi + \mathcal L\psi),\\
  \label{eq:psi-3}
  \psi_{\sym3 H} =\tfrac16 (3\psi - \mathcal L\psi).
\end{gather}
The component $\psi_H$ is characterised by $\mathcal L_A \psi_H = \psi_H$,
for $A=I,J,K$.  On the other hand $\psi_{\sym3H}$~satisfies $\sum_{A=I,J,K}
A_{(12)} \psi_{\sym3H} = -\psi_{\sym3H}$.  Writing $\psi_H = \psi^{(KH)} +
\psi^{(EH)} \in KH + EH$ and $\psi_{\sym3H} = \psi^{(33)} + \psi^{(E 3)}
\in \Lambda_0^3 E \sym3 H+ E \sym3 H$, one computes
\begin{gather*}
  \psi^{(EH)} = - \tfrac1{2n+1} \sum_{A=I,J,K} A \theta^\psi \wedge
  \omega_A,\\
  \psi^{(E3)} = - \tfrac1{2(n-1)} \sum_{A=I,J,K} A \bigl(\theta^\psi_A-
  \theta^\psi\bigr) \wedge \omega_A,
\end{gather*}
where $\theta^\psi_A(X) = A\Lambda_A\psi$ and $\theta^\psi = \frac13
\sum_{A=I,J,K} \theta^\psi_A$.

Let us now describe the $\SP(n)\SP(1)$-components of the intrinsic
torsion~$\xi$, and include characterisations via three-forms.  We will
write $\xi_{33}$, $\xi_{K3}$, $\xi_{E3}$, $\xi_{3H}$, $\xi_{KH}$ and
$\xi_{EH}$ for the components of~$\xi$ corresponding to the modules in the
sum~\eqref{eq:swann-sum}.  We have the following descriptions:
\begin{asparaenum}
  \setdefaultleftmargin{}{4em}{}{}{}{} 
\item $\xi_{33}$ is a tensor characterised by the conditions:
  \begin{enumerate}
  \item $\sum_{A=I,J,K} (\xi_{33})_A A = - \sum_{A=I,J,K} A (\xi_{33})_A =
    - \xi_{33}$,
  \item $\inp\cdot{(\xi_{33})_\cdot\cdot}$ is a skew-symmetric three-form.
  \end{enumerate}
  Or equivalently, $\xi_{33}$ is give by $ \inp Y{(\xi_{33})_XZ}=
  \psi^{(3)}(X,Y,Z)$, where $\psi^{(3)}$ lies in the module $\Lambda_0^3 E
  \sym3 H \subset \Lambda^3 T^*M$.
\item $\xi_{K3}$ is a tensor characterised by the conditions:
  \begin{enumerate}
  \item $\sum_{A=I,J,K} (\xi_{K3})_A A = -\sum_{A=I,J,K} A (\xi_{K3})_A = -
    \xi_{K3}$,
  \item $\sumcic_{XYZ} \inp Y{(\xi_{K3})_X Z} =0$.
  \end{enumerate}
  Or equivalently, $\xi_{K3}$ is expressed by
  \begin{equation*}
    \inp Y{(\xi_{K3})_X Z} = \sum_{A=I,J,K} A_{(23)} \psi^{(K)}_A,
  \end{equation*}
  where $\psi^{(K)}_A$, $A=I,J,K$, are local three-forms in the module $KH$
  such that $\sum_{A=I,J,K} \psi^{(K)}_A = 0$.
\item $\xi_{E3}$ is given by
  \begin{multline*}
    \inp Y{(\xi_{E3})_X Z} \\
    = \tfrac1n\sum_{A=I,J,K} \bigl( nA ( \theta^\xi_A - \theta^\xi ) \wedge
    \omega_A - (n-1) A ( \theta^\xi_A - \theta^\xi) \otimes
    \omega_A\bigr)(X,Y,Z),  
  \end{multline*}
  where $\theta^\xi$ is the one-form defined by
  \begin{equation}
    \label{eq:global-theta}
    \tfrac6n (2n+1)(n-1) \theta^\xi (X) = - \inp{\xi_{e_i} e_i}X = -
    \sum_{A=I,J,K} \inp{A\xi_{e_i} A e_i}X,
  \end{equation}
  and $\theta^\xi_I$, $\theta^\xi_J$, $\theta^\xi_K$ are the local
  one-forms given by
  \begin{equation*}
    \tfrac2n (2n+1) (n-1) \theta^\xi_A (X) = - \inp{A \xi_{e_i} A e_i}X.
  \end{equation*}
  Note that $3\theta^\xi = \theta^\xi_I + \theta^\xi_J + \theta^\xi_K$.
\item $\xi_{3H}$ is a tensor characterised by the conditions:
  \begin{enumerate}
  \item $(\xi_{3H})_A A - A (\xi_{3H})_A -A \xi_{3H} A = \xi_{3H}$, for
    $A=I,J,K$,
  \item $\sumcic_{X,Y,Z} \inp Y{(\xi_{3H})_X Z} = 0$.
  \end{enumerate}
  Or equivalently, $\xi_{3H}$ is expressed by
  \begin{equation*}
    \inp Y{(\xi_{3H})_X Z} = \sum_{A=I,J,K} A_{(23)} \psi^{(3)}_A,
  \end{equation*}
  where $\psi^{(3)}_A$, $A=I,J,K$ are local
  three-forms such that
  \begin{enumerate}
    \setcounter{enumii}{15} 
  \item $\psi^{(3)}_A$ is
    in $\Lambda_0^3 E \sym3 H$,
  \item $\psi^{(3)}_A$ is of type $\{2,1\}$ with respect to the almost
    complex structure~$A$, i.e., $\mathcal L_A \psi^{(3)}_A=\psi^{(3)}_A$,
    $A=I,J,K$, and
  \item $\sum_{A=I,J,K} \psi^{(3)}_A = 0$.
  \end{enumerate}
  One may check that one of these three-forms is sufficient to determine
  the others.  Indeed
  \begin{equation*}
    \psi^{(3)}_J = -\tfrac12(3+\mathcal L_J)\psi^{(3)}_I
    \quad\text{and}\quad
    \psi^{(3)}_K = -\tfrac12(3+\mathcal L_K)\psi^{(3)}_I.
  \end{equation*}
\item $\xi_{KH}$ is a tensor characterised by the conditions:
  \begin{enumerate}
  \item\label{item:KH-1} $ (\xi_{KH})_A A - A (\xi_{KH})_A - A \xi_{KH} A =
    \xi_{KH}$, for $A=I,J,K$;
  \item\label{item:KH-2} there exists a skew-symmetric three-form
    $\psi^{(K)}$ such that
    \begin{equation*}
      \inp Y{(\xi_{KH})_X Z} =\bigl(3\psi^{(K)} - \sum_{A=I,J,K} A_{(23)}
      \psi^{(K)} \bigr)(X,Y,Z);
    \end{equation*}
  \item\label{item:KH-3} $\sum_{i=1}^{4n} (\xi_{KH})_{e_i} e_i =0$.
  \end{enumerate}
  Note that conditions \itref{item:KH-1} and \itref{item:KH-3} can be
  replaced by saying that $\psi^{(K)}$~is in~$KH$, i.e., for each
  $A=I,J,K$, the form~$\psi^{(K)}$ is of type $\{2,1\}_A$ and satisfies
  $\Lambda_A\psi^{(K)} = 0$.
\item $\xi_{EH}$ is given by
  \begin{equation*}
    \begin{split}
      \inp Y{(\xi_{EH})_X Z} &= 3 e_i \otimes e_i \wedge \theta^\xi(X,Y,Z)
      \eqbreak - \sum_{A=I,J,K} \bigl( e_i \otimes A e_i \wedge A
        \theta^\xi + \tfrac2{n} A \theta^\xi \otimes \omega_A \bigr)
      (X,Y,Z), 
    \end{split}
  \end{equation*}
  where $\theta^\xi$ is the global one-form defined by~\eqref{eq:global-theta}.
\item The part $\xi_{\sym3H} = \xi_{33} + \xi_{K3} + \xi_{E3} $ of $\xi$ in
  $( \Lambda_0^2 E + K + E) \sym3H$ is characterised by the condition
  \begin{equation*}
    \sum_{A=I,J,K} (\xi_{\sym3 H })_A A = - \sum_{A=I,J,K} A
    (\xi_{\sym3H})_A = - \xi_{\sym3H} . 
  \end{equation*}
\item The part $\xi_H = \xi_{3H} + \xi_{KH} + \xi_{EH} $ of $\xi$ in $(
  \Lambda_0^2 E + K + E) H$ is characterised by the condition
  \begin{equation*}
    (\xi_{ H })_A A - A (\xi_{H})_A - A (\xi_{H}) A = \xi_H,
  \end{equation*}
  for $A=I,J,K$.
\end{asparaenum}

\section{Use of exterior derivatives}
\label{sec:exterior}

Here we will find out how the intrinsic torsion $\xi$ can be determined by
means of the exterior derivatives of the K\"ahler forms $d \omega_I$, $d
\omega_J$ and $d \omega_K$.

In \cite{Cabrera-S:aHqg} it was shown that the covariant derivatives
$\nabla \omega_I$, $\nabla \omega_J$ and $\nabla \omega_K$ are given by
\begin{equation}
  \label{eq:lambda-alpha}
  \LC \omega_I = \lambda_K \otimes \omega_J - \lambda_J \otimes \omega_K +
  J_{(2)} \alpha_K - K_{(2)} \alpha_J, \detc 
\end{equation}
where $\lambda_A$ are given by equation~\eqref{eq:lambda-I} and $\alpha_I,
\alpha_J , \alpha_K \in T^*M \otimes \Lambda^2_0 E \subset T^*M \otimes
\sym2 T^*M$ are defined by
\begin{equation}
  \label{eq:alpha-I}
  \begin{split}
    \alpha_I
    &:= - \lambda_I \otimes g + \tfrac12 ( J_{(2)} -J_{(3)}) \LC \omega_K
    \\ 
    &= - \lambda_I \otimes g + \tfrac12 ( K_{(3)} -K_{(2)}) \LC \omega_J,
    \detc
  \end{split}
\end{equation}
We may rewrite the intrinsic torsion~$\xi$ from
Proposition~\ref{prop:torsion} using equation~\eqref{eq:lambda-alpha}
giving
\begin{equation}
  \label{eq:xi-alpha}
  \xi_XY = - \tfrac12 \bigl((\alpha_I)_X IY) + (\alpha_J)_X JY +
  (\alpha_K)_X KY \bigr),  
\end{equation}
where $\left( \alpha_A \right)$ is given by $\inp Y{\left( \alpha_A
\right)_X Z} = \alpha_A ( X ; Y , Z)$.  Thus the intrinsic torsion $\xi$
only depends on the $\alpha_A$'s; the $\lambda_A$'s have no influence.

Note that the dimension of the space of possible triples $(\alpha_I,
\alpha_J, \alpha_K)$ coincides with the dimension $\dim T^*M \otimes
\Lambda_0^2 E \sym2 H = 12n (2n+1)(n-1)= 3 \dim T^*M \otimes \Lambda^2_0 E
$ of the space of possible intrinsic torsion tensors.

As $EH \otimes \Lambda_0^2 E = \Lambda_0^3 EH + KH + EH $, one may
decompose $\alpha_I$ into three $\SP(n)\SP(1)$-components
\begin{equation*}
  \alpha_I = \alpha_I^{(3)} + \alpha_I^{(K)} + \alpha_I^{(E)} \in
  \Lambda_0^3 EH + KH + EH.
\end{equation*}
If $\dim M=8$, the module $\Lambda^3_0E$ is trivial and the corresponding
component~$\alpha_I^{(3)}$ is not present.  The component $\alpha_I^{(E)}$
is determined from a one-form $\eta_I$ which is defined by
\begin{equation}
  \label{eq:eta-I}
  \eta_I (X) = \alpha_I(e_i,e_i,X).
\end{equation}
Furthermore, the components of the $\alpha_A$'s can be used to characterise
classes of almost quaternion-Hermitian manifolds.

\begin{proposition}[Cabrera \& Swann~\cite{Cabrera-S:aHqg}]
  \label{prop:alpha}
  Let $M$ be an almost quaternion-Hermitian manifold with intrinsic
  torsion~$\xi$.  If $I,J,K$ is an adapted basis of $\mathcal G$, then for
  $V=3, K, E$,
  \begin{enumerate}
  \item each component $\xi_{VH}$ is linearly determined by
    $I_{(1)}\alpha^{(V)}_I + J_{(1)}\alpha^{(V)}_J +
    K_{(1)}\alpha^{(V)}_K$,
  \item each component $\xi_{V3}$ is linearly determined by
    $I_{(1)}\alpha^{(V)}_I - J_{(1)}\alpha^{(V)}_J$ and
    $J_{(1)}\alpha^{(V)}_J - K_{(1)}\alpha^{(V)}_K$.
  \end{enumerate}
\end{proposition}

Observing that $A_{(1)}\alpha^{(E)}_A=A\alpha^{(E)}_A$ is linearly
determined by the one-form $A \eta_A$, we have the following result.

\begin{corollary}
  \label{cor:E}
  Under the same conditions as Proposition~\ref{prop:alpha}, we have:
  \begin{enumerate}
  \item $\xi_{EH}$ is linearly determined by $I\eta_I+J\eta_J+K\eta_K$,
  \item $\xi_{E3}$ is linearly determined by $I\eta_I -J\eta_J$ and
    $J\eta_J-K\eta_K$.
  \end{enumerate}
\end{corollary}

We now proceed to express the covariant derivatives $\LC_X \omega_I$,
$\LC_X \omega_J$ and $\LC_X \omega_K$ in terms of $d \omega_I$, $d
\omega_J$ and~$d \omega_K$.  We use a relation between these covariant
derivatives found in~\cite{Fernandez-Moreiras:symmetry,Cabrera:aqh}, which
may be symmetrically expressed by
\begin{equation}
  \label{eq:sym-nabla}
  \left(\LC_X \omega_I \right)(JY,KZ) + \left(\LC_X \omega_J \right)(KY,IZ)
  + \left(\LC_X \omega_K \right)(IY,JZ) = 0
\end{equation}
and the following identity given by Gray~\cite{Gray:minimal}
\begin{equation}
  \label{eq:Gray}
  2 \LC\omega_I = d\omega_I - I_{(23)} d\omega_I - I_{(3)} N_I,
\end{equation}
where $N_I (X,Y,Z) = \inp X{N_I(Y,Z)}$ and the $(1,2)$-tensor $N_I$ is the
Nijenhuis tensor for~$I$, i.e., $N_I(X,Y) = [X,Y] + I[IX,Y] + I[X,IY] -
[IX,IY]$.

Under the action of~$\Un(2n)_I$ the space of three-forms decomposes in to
irreducible modules as
\begin{equation*}
  \Lambda^3 T^*M = \Lambda^{\{3,0\}}_IT^*M + \Lambda^{\{2,1\}}_{0,I}T^*M +
  \Lambda^{\{1,0\}}_IT^*M\wedge\omega_I = \Wc{1+3+4,I}, 
\end{equation*}
where $\Wc{iI}$ are isomorphic to the Gray-Hervella modules
described in~\cite{Gray-H:16} and the subscript~$I$ indicates the almost
complex structure considered.  Note that a three-form~$\psi$ lies in
$\Wc{3+4,I} = \Lambda^{\{2,1\}}_IT^*M \subset \Lambda^3 T^*M$ if and only if
\begin{equation*}
  (I_{(12)} + I_{(13)} + I_{(23)} ) \psi = \psi,
  \quad\text{i.e.,}\ \mathcal L_I \psi = \psi.
\end{equation*}

\begin{proposition}
  \label{prop:naddd}
  For an adapted basis $I,J,K$ the exterior derivatives $d\omega_I$, etc.,
  determine
  \begin{asparaenum}
  \item\label{item:naddd-1} the covariant derivative $\LC \omega_I$ by
    \begin{align}
      \label{naddd3}
      2 \LC \omega_I &= (1 - I_{(23)}) d \omega_I + (I_{(2)} + I_{(3)})
      J_{(1)} d \omega_J \eqbreak\nonumber\qquad - (1 - I_{(23)} ) J_{(1)}
      d \omega_K
      \\
      \label{naddd4}
      &= (1 - I_{(23)}) d \omega_I + (I_{(2)}+I_{(3)})K_{(1)} d \omega_K
      \eqbreak\nonumber\qquad + (1 - I_{(23)} ) K_{(1)} d \omega_J,
    \end{align}
  \item\label{item:naddd-2} the Nijenhuis $(0,3)$-tensor $N_I$ by
    \begin{equation}
      \label{eq:naddd2}
      \begin{split}
        2 N_I & = ( I_{(12)} + I_{(13)} + I_{(23)} - 1) J_{(23)} (J d
        \omega_J - K d \omega_K)\\
        & = ( 1- I_{(12)} ) ( K_{(23)} - J_{(23)}) (J d \omega_J - K d
        \omega_K),
      \end{split}
    \end{equation}
  \item\label{item:naddd-3} the one-form $I\lambda_I$ of
    equation~\eqref{eq:lambda-I} by
    \begin{align}
      \label{eq:ladd1}
      2n I\lambda_I & = I\Lambda_K d \omega_J +J\Lambda_K d \omega_I + J
      \Lambda_I d \omega_K\\  
      \label{eq:ladd2}
      & =  I\Lambda_K d \omega_J + \Lambda_I d \omega_I - \Lambda_K d \omega_K, 
    \end{align}
    and moreover
  \item\label{item:naddd-4} $J d \omega_J + K d \omega_K \in \Wc{3+4,I}$, i.e.,
    \begin{gather}
      \label{eq:naddd4r}
      (I_{(12)} + I_{(13)} + I_{(23)} ) (J d \omega_J + K d \omega_K) = J d
      \omega_J + K d \omega_K.
    \end{gather}
  \end{asparaenum}
  The corresponding expressions with respect to $J$ and $K$ are obtained by
  cyclic permutations of $I,J,K$.
\end{proposition}

\begin{proof}
  Equation~\eqref{naddd3} is derived from its right-hand side, taking into
  account that $ d \omega_A (X,Y,Z) = \sumcic_{XYZ} \left( \LC_X \omega_A
  \right) (Y,Z)$, $A=I,J,K$, and making repeated use of
  equation~\eqref{eq:sym-nabla}.  The proof for equation~\eqref{naddd4} is
  similar.  Now, \itref{item:naddd-2}, \itref{item:naddd-3} and
  \itref{item:naddd-4} are immediate consequences of~\eqref{naddd3},
  \eqref{naddd4} and Gray's identity~\eqref{eq:Gray}.
\end{proof}

The expressions for the intrinsic $\SP(n)\SP(1)$-torsion $\xi$ given in
next result are consequences of the last proposition and
equation~\eqref{eq:aqH-torsion}.

\begin{proposition}
  The intrinsic $\SP(n)\SP(1)$-torsion $\xi$ is determined by the exterior
  derivatives $d\omega_I$, etc., by
  \begin{multline*}
      \xi_X Y = \tfrac1{4n} \sumcic_{IJK} (\Lambda_K d
        \omega_J - I\Lambda_I d \omega_I + I\Lambda_K d \omega_K)(X)IY \\
      +\tfrac18 \sumcic_{IJK}
      \Bigl(\bigl(I_{(2)}+I_{(3)}+(J_{(12)}+J_{(13)}+K_{(23)}-1)I_{(1)}\bigr)
      d\omega_I\Bigr)(X,Y,e_i)e_i.  
  \end{multline*}
\end{proposition}

\begin{proof}
  One computes first
  \begin{equation*}
    \begin{split}
      \xi_X Y &= \tfrac1{4n} \sumcic_{IJK} \left\{ \inp{X \hook d
        \omega_J}{\omega_K} + \inp{I X \hook d \omega_I}{\omega_I} -
        \inp{IX \hook d \omega_K}{\omega_K}
      \right\} IY \\
      &\ - \tfrac18 \sumcic_{IJK} \bigl\{ d \omega_I (X,Y,Ie_i) + d
      \omega_I (X,IY,e_i) - d \omega_J (JX,Y,e_i)\eqbreak + d \omega_J
      (JX,IY,Ie_i)+ d \omega_K (JX,Y,Ie_i)e_i + d \omega_K (JX,IY,e_i)
      \bigr\} e_i;
    \end{split}
  \end{equation*}
  and then takes advantage of the second cyclic sum to rearrange terms.
\end{proof}

\section{A minimal description}
\label{sec:beta}

Motivated by Proposition~\ref{prop:naddd}\itref{item:naddd-4}, let us
introduce the three-forms 
\begin{equation*}
  \beta_I = J d \omega_J + K d \omega_K, \detc
\end{equation*}
These determine the exterior derivatives $d\omega_I$, $d\omega_J$,
$d\omega_K$ as follows
\begin{equation*}
  2 d \omega_I = I (\beta_I -  \beta_J -  \beta_K), \detc.
\end{equation*}
and $\beta_A \in \Wc{3+4,A}=\Lambda^{\{2,1\}}_AT^*M$, so the dimension of
the space of possible exterior derivatives $d\omega_I$, $d\omega_J$,
$d\omega_K$ is at most $3 \dim (\Wc{3+4}) = 12 n^2 (2n-1)$
\cite{Gray-H:16}.  On the other hand, equation~\eqref{eq:lambda-alpha}
implies that the dimension of the space of covariant derivatives $\nabla
\omega_I$, $\nabla \omega_J$, $\nabla \omega_K$ is determined by the
possible triples of $\lambda$'s and $\alpha$'s.  This dimension is $12 n +
12n (2n+1)(n-1) = 12 n^2 (n-1)$, which coincides with the above one
computed for the~$\beta$'s.  Therefore, algebraically, the three-forms
$\beta_I$, $\beta_J$, $\beta_K$ are independent.

We will now show how, components of, the $\beta$'s determine the intrinsic
torsion.  Consider the action of the group $\SP(n)\Un(1)_I$, which is the
intersection of $\Un(2n)_I$ with~$\SP(n)\SP(1)$, on the module $\Wc{3+4,I}
= \Lambda^{\{2,1\}}_IT^*M \subset \Lambda^3 T^*M$.  It was shown
in~\cite{Cabrera-S:aHqg} that $\Wc3 \otimes \mathbb C = (\Lambda_0^3 E + K
+ E) (L_I + \overline{L_I}) $ and $\Wc4 \otimes \mathbb C = E
(L_I+\overline{L_I})$, where we write $L_I$ for the standard representation
of $\Un(1)_I$ on~$\mathbb C$.  Since $\Lambda_0^3 E$, $K$, $E$ are
representations of quaternionic type and $L_I$~is of complex type, the
tensor products $\Lambda^3_0EL$, etc., in the above decompositions are all
of quaternionic type.  The underlying real modules $[V]_{\mathbb R}$
obtained by regarding the modules $V$ as real vector spaces give real
representations of~$\SP(n) \Un(1)$ and
\begin{equation*}
  \Wc{3I} = [\Lambda_0^3 EL_I]_{\mathbb R} + [KL_I]_{\mathbb R} + [EL_I]_{\mathbb
  R3}, \qquad \Wc{4I} = [EL_I]_{\mathbb R4}.
\end{equation*}
Using these decompositions, the tensor~$\beta_I$ splits into four
components
\begin{equation}
  \label{eq:beta}
  \beta_I = \beta_I^{(3)} + \beta_I^{(K)} + \beta_{3I}^{(E)}
  + \beta_{4I}^{},
\end{equation}
with one-form parts
\begin{gather}
  \label{eq:beta-E3}
  \beta_{3I}^{(E)} =  - \tfrac12 J \nu^I_3 \wedge \omega_J - \tfrac12 K
  \nu^I_3 \wedge \omega_K + \tfrac1{2n-1} I \nu^I_3 \wedge \omega_I,\\  
  \beta_{4I} =  - \tfrac1{2n-1} I \nu^I_4 \wedge \omega_I,
  \label{eq:beta-E4}
\end{gather}
where $\nu^I_3$ and $\nu^I_4$ are one-forms, which we will now specify.  We
have
\begin{equation}
  \label{eq:nu-I4}
  \nu^I_4 = I\Lambda_I\beta_I.
\end{equation}
A computation gives the following formula determining $\nu^I_3$
from~$\beta_I$:
\begin{equation}
  \label{eq:beta-theta}
  J\Lambda_J\beta_I = K\Lambda_K\beta_I =
  \tfrac1{2n-1}\bigl(\nu^I_4 + (2n+1)(n-1) \nu^I_3\bigr).
\end{equation}
Here the first equality in~\eqref{eq:beta-theta} is equivalent to
\begin{equation}
  \label{eq:ijk}
  I\Lambda_K d\omega_J + I\Lambda_J d \omega_K = - \Lambda_J d\omega_J +
  \Lambda_K d\omega_K
\end{equation}
which is an immediate consequence of equations \eqref{eq:ladd1}
and~\eqref{eq:ladd2} of Proposition~\ref{prop:naddd} and the fact that
$\beta_A\in \Wc{3+4,A}$.  We may now find the other components of~$\beta_I$
via \eqref{eq:psi-H} and~\eqref{eq:psi-3}:
\begin{equation}
  \label{eq:beta-K3}
  \beta_I^{(3)} = \tfrac16 ( 2 - \mathcal
  L_J - \mathcal L_K) \beta_I^{(3+K)} , \quad \beta_I^{(K)} = \tfrac16
  ( 4 + \mathcal L_J + \mathcal L_K) \beta_I^{(3+K)},
\end{equation}
where $\beta_I^{(3+K)} = \beta_I^{} - \beta_{3I}^{(E)} - \beta_{4I}^{}$.

\begin{remark}
  The expressions for the one-form parts are a little simpler in dimension
  four, i.e., $n=1$.  Recall that the Lee form of~$\omega_I$ is $Id^*
  \omega_I = - \Lambda_I d\omega_I$, where $d^*$ is the co-derivative.  For
  $n=1$, we have $\Wc3=\{0\}$, so $\beta_I \in \Wc{4I}$, $I\Lambda_I\beta_I
  = J\Lambda_J\beta_I = K\Lambda_K\beta_I = \nu^I_4$, and
  \begin{equation*}
    K\Lambda_J d\omega_I = - J\Lambda_K d\omega_I = - \Lambda_I d\omega_I =
    Id^*\omega_I, \detc  
  \end{equation*}
\end{remark}

\begin{remark}
  In order to apply Proposition~\ref{prop:alpha} to classify almost
  quaternion-Hermitian manifolds, the tensors $\alpha_A^{(3)}$ and
  $\alpha_A^{(K)}$ can be computed from the triples $\beta_I^{(3)}$,
  $\beta_J^{(3)}$, $\beta_K^{(3)}$ and $\beta_I^{(K)}$, $\beta_J^{(K)}$,
  $\beta_K^{(K)}$ respectively.  In fact, we would begin with
  equations~\eqref{eq:alpha-I} which define~$\alpha_A$ and then use
  Proposition~\ref{prop:naddd}.
\end{remark}

To analyse the $E\sym3H$ and $EH$ components of the intrinsic torsion
$\xi$, we wish to apply Corollary~\ref{cor:E} which requires knowledge of
the one-forms~$\eta_I$.  Let us see how these are determined by the
$\beta_I$'s. Equation~\eqref{eq:lambda-alpha} gives
\begin{equation*}
  \begin{split}
    \beta_I &= - J \lambda_I \wedge \omega_K + J \lambda_K \wedge
    \omega_I + J \Alt (K_{(2)} \alpha_I) - J \Alt (I_{(2)} \alpha_K) \eqbreak
    - K \lambda_J \wedge \omega_I + K \lambda_I \wedge \omega_J + K \Alt
    (I_{(2)} \alpha_J) - K \Alt (J_{(2)} \alpha_I),
  \end{split}
\end{equation*}
where $\Alt(\phi)(X,Y,Z) = \sumcic_{XYZ} \phi(X,Y,Z)$, for $\phi \in T^*M
\otimes \Lambda^2 T^*M$.  This combined with~\eqref{eq:ladd2} leads to
\begin{gather*}
  \begin{split}
    4n I\eta_I &= 2(n-1)J\Lambda_J\beta_I +
    I\Lambda_I((n-1)\beta_I+\beta_J+\beta_K) \eqbreak - n J\Lambda_J\beta_J
    - n K\Lambda_K\beta_K,
  \end{split}
  \\
  4n I\lambda_I = 2J\Lambda_J\beta_I + I\Lambda_I(\beta_I-\beta_J-\beta_K),
  \detc
\end{gather*}
Note that the right-hand sides of these equations are linear combinations
of~$\nu^A_3$ and $ \nu^A_4$, $A=I,J,K$, so
\begin{align*}
  I \eta_I & =  \tfrac{(2n+1)(n-1)}{4n (2n-1)} \, \bigl( ( 2 (n-1) \nu^I_3
  + \nu^J_3 + \nu^K_3) + ( \nu^I_4 - \nu^J_4 - \nu^K_4)\bigr),   
  \\
  I \lambda_I & =  \tfrac{(2n+1)(n-1)}{4n(2n-1)} \, ( 2 \nu^I_3 - \nu^J_3 -
  \nu^K_3) + \tfrac1{4n(2n-1)} ((2n+1) \nu^I_4 - \nu^J_4 - \nu^K_4), 
\end{align*}
etc.

The next proposition shows clearly the r\^oles played by the three-forms
$\beta_I$, $\beta_J$ and $\beta_K$ in determining the components of the
intrinsic torsion~$\xi$.  This provides a practical way to compute $\xi$
using the tools of the exterior algebra.

\begin{proposition}
  \label{prop:torsion-beta}
  For an almost quaternion-Hermitian $4n$-manifold, $n>1$, we have:
  \begin{enumerate}
  \item\label{item:tor-beta-1} The three-form $\psi^{(3)}$, which
    determines $\xi_{33}$, is given by
    \begin{equation}
      \label{beta33}
      \psi^{(3)} = \tfrac1{12} (\beta^{(3)}_I + \beta^{(3)}_J + \beta^{(3)}_K).
    \end{equation}
  \item\label{item:tor-beta-2} The local three-forms $\psi^{(3)}_I$,
    $\psi^{(3)}_J$, $\psi^{(3)}_K$ each of which determine $\xi_{3H}$, are
    given by
    \begin{equation}
      \label{eq:beta-3H}
      \psi^{(3)}_A
      = - \tfrac18 \beta^{(3)}_A + \tfrac1{48}(3+\mathcal L_A) \sum_{B=I,J,K}
      \beta^{(3)}_B .
    \end{equation}
  \item\label{item:tor-beta-3} The three-form $\psi^{(K)}$, which
    determines $\xi_{KH}$, is given by
    \begin{equation}
      \label{eq:beta-KH}
      \psi^{(K)} = - \tfrac1{48} (\beta^{(K)}_I + \beta^{(K)}_J +
      \beta^{(K)}_K ). 
    \end{equation}
  \item\label{item:tor-beta-4} The local three-forms $\psi^{(K)}_I$,
    $\psi^{(K)}_J$, $\psi^{(K)}_K$ which determine~$\xi_{K3}$, are given by
    \begin{equation}
      \label{eq:beta-K3-2}
      \psi^{(K)}_A = -\tfrac12 \beta^{(K)}_A + \tfrac16\sum_{B=I,J,K}
      \beta^{(K)}_B.
    \end{equation}
  \item\label{item:tor-beta-5} The one-form $\theta^\xi$, which determines
    $\xi_{EH}$, is given by
    \begin{equation}
      \label{eq:theta-xi}
      \theta^\xi = \tfrac n{24(2n-1)} \sum_{A=I,J,K} ( \nu^A_3 -
      2 A \lambda_A), 
    \end{equation}
  \item\label{item:tor-beta-6} The local three-forms $\theta^\xi_I$,
    $\theta^\xi_J$, $\theta^\xi_K$ whose differences $\theta_A^\xi -
    \theta^\xi$ determine $\xi_{E3}$, are given by
    \begin{equation}
      \label{eq:theta-IJK}
      \theta^\xi_A = \tfrac n{4(n+1)}\Bigl((\nu^A_3 - 2 A
      \lambda_A) - \tfrac{n-1}{2(2n-1)} \sum_{B=I,J,K} (\nu^B_3 -
      2 B \lambda_B)\Bigr).
    \end{equation}
  \end{enumerate}
\end{proposition}

\begin{proof}
  For the covariant derivative of the local K\"ahler forms~$\omega_I$, we
  have
  \begin{equation}
    \label{eq:LC-xi}
    \begin{split}
      \LC_X \omega_I(Y,Z) &= \lambda_K(X) \omega_J(Y,Z) - \lambda_J(X)
      \omega_K(Y,Z) \eqbreak - \inp Y{\xi_XIZ} - \inp{IY}{\xi_XZ},
    \end{split}
  \end{equation}
  from which one derives
  \begin{equation*}
    I d \omega_I = - I \lambda_K \wedge
    \omega_J + I\lambda_J \wedge \omega_K - \sumcic_{XYZ}
    \bigl(\inp Y{\xi_{IX}IZ} + \inp{IY}{\xi_{IX}Z}\bigr)
  \end{equation*}
  and
  \begin{equation}
    \label{eq:beta-xi}
    \begin{split}
      \beta_I & = K\lambda_I \wedge \omega_J - J\lambda_I \wedge
      \omega_K -I \lambda_I^+ \wedge \omega_I \eqbreak
      - \sumcic_{XYZ} \bigl(\inp Y{\xi_{JX}JZ} + \inp{JY}{\xi_{JX}Z} + \inp
      Y{\xi_{KX}KZ} + \inp{KY}{\xi_{KX}Z}\bigr), 
    \end{split}
  \end{equation}
  where $\lambda^+_I = J \lambda_J + K \lambda_K$.

  For parts \itref{item:tor-beta-1} and \itref{item:tor-beta-2},
  equation~\eqref{eq:beta-xi} gives
  \begin{equation*}
    \begin{split}
      -\beta^{(3)}_I (X,Y,Z) = \sumcic_{XYZ} &\inp Y{(\xi_{33} +
      \xi_{3H})_{JX}JZ} + \inp{JY}{(\xi_{33}+\xi_{3H})_{JX}Z} \\
      &+ \inp Y{(\xi_{33}+\xi_{3H})_{KX} KZ} +
      \inp{KY}{(\xi_{33}+\xi_{3H})_{KX}Z}.  
    \end{split}
  \end{equation*}
  We now get
  \begin{gather*}
    \beta^{(3)}_A = 6 \psi^{(3)}+ 2 \mathcal L_A \psi^{(3)} - 8
    \psi^{(3)}_A
  \end{gather*}
  which leads to equations \eqref{beta33} and~\eqref{eq:beta-3H} as required.
  
  For parts \itref{item:tor-beta-3} and \itref{item:tor-beta-4}, we
  use~\eqref{eq:beta-xi} to get
  \begin{gather*}
    \beta^{(K)}_A =- 16 \psi^{(K)} - 2 \psi^{(K)}_A,
  \end{gather*}
  which gives equations \eqref{eq:beta-KH} and~\eqref{eq:beta-K3-2}.
  
  Finally, for parts \itref{item:tor-beta-5} and \itref{item:tor-beta-6},
  use \eqref{eq:beta-xi} to find
  \begin{equation*}
    \label{eq:beta-E}
    \begin{split}
      n\beta^{(E)}_I &=-I \bigl(n \lambda_I^+ + 6(n-1) \theta^\xi -
      2(3n+1) \theta^\xi_I \bigr) \wedge \omega_I \eqbreak
        - J \bigl(n I\lambda_I + 6(n-1) \theta^\xi
        + 2(n+1) \theta^\xi_I \bigr) \wedge
        \omega_J \eqbreak
        - K \bigl(n I\lambda_I + 6(n-1) \theta^\xi +
        2(n+1) \theta^\xi_I \bigr) \wedge \omega_K, \detc
    \end{split}
  \end{equation*}
  Using equations \eqref{eq:beta-E3}, \eqref{eq:beta-E4}
  and~\eqref{eq:beta-theta}, this gives
  \begin{align}
    \label{eq:beta-3A}
    \nu^A_3 & =  2 A\lambda_A + \tfrac4{n} \bigl(3(n-1)\theta^\xi +
    (n+1)\theta^\xi_A\bigr), \\ 
    \label{eq:beta-4A}
    \nu^A_4 & = (2n-1) \lambda_A^+ + 2 A\lambda_A +
    \tfrac{6(n-1)(2n+1)}{n} (\theta^\xi - \theta^\xi_A), 
  \end{align}
  for $A=I,J,K$.  As $3 \theta^\xi = \theta^\xi_I + \theta^\xi_J +
  \theta^\xi_K$, equations \eqref{eq:theta-xi} and~\eqref{eq:theta-IJK}
  follow.
\end{proof}

We may now quickly record what happens under conformal changes of metric.

\begin{proposition}
  \label{prop:conformal}
  On a almost quaternion-Hermitian $4n$-manifold, if we consider a
  conformal change of metric $\inp{\cdot}{\cdot}^o = e^{2\sigma}
  \inp{\cdot}{\cdot}$, with $\sigma\in C^\infty(M)$, then
  \begin{gather*}
    \omega_A^o = e^{2\sigma} \omega_A, \quad d \omega_A^o = e^{2\sigma}
    \left( 2 d \sigma \wedge \omega_A + d \omega_A
    \right),\\
    A d^* \omega_A^o = A d^* \omega_A - 2 (2n-1) d\sigma, \quad
    A \lambda_A^o = A \lambda_A - \tfrac1{n} d\sigma, \\
    \beta_I^o = e^{2 \sigma} \left( \beta_I + 2 J d \sigma \wedge \omega_J
      + 2 K d \sigma \wedge \omega_K \right), \detc, \\
    {\nu^A_3}^o = \nu^A_3 - 4 d \sigma, \quad
    {\nu^A_4}^o = \nu^A_4 - 4 d \sigma, \quad
    {\theta^\xi}^o = \theta^\xi - \tfrac14 d \sigma, \quad \theta^{\xi^o}_A
    = \theta^\xi_A - \tfrac14 d \sigma.
  \end{gather*}
  In particular, the only component of the intrinsic torsion that changes
  is~$\xi^o_{EH}$.
\end{proposition}

\begin{proof}
  The identities follow from the definitions of each tensor involved.  For
  the intrinsic torsion, use these identities,
  Proposition~\ref{prop:torsion-beta} and the descriptions of the
  components of~$\xi$ given at the end of~\S\ref{sec:intrinsic}.
\end{proof}

\section{HyperK\"ahler manifolds with torsion}
\label{sec:HKT}

In this section we will see some consequences of
Proposition~\ref{prop:naddd} in HKT-geometry.  This geometry arises on the
target space of a $N=2$ super-symmetric $(4,0)$ $\sigma$-models with
Wess-Zumino term.

\begin{definition}[Howe \& Papadopoulos~\cite{Howe-P:twistor-kaehler}]
  \label{def:HKT}
  An almost hyperHermitian manifold $(M,I,J,K,g=\inp\cdot\cdot)$ is an
  \emph{HKT-manifold} (hyperK\"ahler with torsion), if the following
  conditions are satisfied:
  \begin{enumerate}
  \item the almost complex structures $I,J,K$ are integrable;
  \item $M$ admits a linear connection $\Nhkt= \LC + \frac12 T$, such that
    \begin{enumerate}
    \item $\Nhkt I= \Nhkt J = \Nhkt K=0$, and
    \item $\Nhkt g=0$;
    \end{enumerate}
  \item the $(0,3)$-tensor field, also denoted by $T$, defined by
    $T(X,Y,Z)= \inp X{T(Y,Z)}$ is a skew-symmetric three-form.
  \end{enumerate}
\end{definition}

A result of Grancharov \& Poon \cite{Grantcharov-P:HKT} says that an almost
hyperHermitian manifold~$M$ is HKT if and only if $(I,J,K, \inp\cdot\cdot)$
is hyperHermitian (i.e., $N_I=N_J=N_K=0$) and $Id \omega_I= Jd \omega_J = K
d \omega_K$.  We now give the following improvement of this result, showing
that the integrability assumption is redundant.

\begin{proposition}
  Let $(M,I,J,K,g)$ be an almost hyperHermitian manifold.  Then the
  following conditions are equivalent:
  \begin{enumerate}
  \item $M$ is an HKT-manifold;
  \item $Id\omega_I = Jd\omega_J = Kd\omega_K$;
  \item $\beta_I = \beta_J = \beta_K$.
  \end{enumerate}
\end{proposition}

\begin{proof}
  If $M$ is a $HKT$-manifold, we have a connection $\Nhkt = \LC + \frac12T$
  satisfying the conditions given in Definition~\ref{def:HKT}.  The
  integrability condition gives $N_I=0=N_J=N_K$ and $\LC\omega_A \in
  \Wc{3+4,A}$.  Now, using equation~\eqref{eq:Gray}, we obtain $T =
  Id\omega_I= Jd\omega_J = Kd\omega_K = \tfrac12 \beta_I \in \Wc{3+4}$.

  Conversely, suppose $Id\omega_I= Jd\omega_J = Kd\omega_K$.
  Proposition~\ref{prop:naddd}\itref{item:naddd-2} gives $N_I = N_J = N_K =
  0$.  The connection $\Nhkt = \LC + \tfrac12 T$, where $\inp X{T(Y,Z)} =
  Id\omega_I(X,Y,Z)$ now satisfies the HKT conditions.
\end{proof}

Grantcharov \& Poon \cite{Grantcharov-P:HKT} give a second characterisation
of HKT manifolds in terms of the complex geometry of~$I$.  Let us
define as usual the operators $\partial_A$ and $\bar\partial_A$
acting on a $p$-form~$\psi$ by
\begin{equation*}
  \partial_A \psi = \tfrac12 \left( d +  (-1)^p \, i A d A \right) \psi,
  \qquad
  \bar\partial_A \psi = \tfrac12 \left( d - (-1)^p \, i A d A \right) \psi.
\end{equation*}
Assuming integrability of $I$, $J$ and~$K$, Grantcharov \& Poon
show that $M$ is HKT if and only if the $(2,0)$-form $\omega_J+i\omega_K$ is $\partial_I$-closed.  Once again we may weaken
the integrability requirements.

\begin{proposition}
  \label{prop:HKT}
  Let $(M,I,J,K,g)$ be an almost hyperHermitian manifold.  Then the
  following conditions are equivalent:
  \begin{enumerate}
  \item $M$ is an HKT-manifold;
  \item $Jd\omega_J = Kd\omega_K$ and $N_J=0$;
  \item $\partial_I(\omega_J + i\omega_K)=0$ and $N_J=0$;
  \item $\bar\partial_I(\omega_J - i\omega_K)=0$ and $N_J=0$.
  \end{enumerate}
\end{proposition}

\begin{proof}
  It is easy to see that the three conditions $ Jd\omega_J = Kd\omega_K$,
  $\partial_I (\omega_J + i\omega_K) = 0$ and $\bar\partial_I (\omega_J -
  i\omega_K) = 0$, are equivalent.  Moreover, by
  Proposition~\ref{prop:naddd}\itref{item:naddd-2}, the condition
  $Jd\omega_J = K d\omega_K$ implies $N_I=0$.  Now, the integrability of
  $I$ and $J$ imply that $K$ is integrable (see \cite{Obata:connection} or
  the newer proof~\cite{Cabrera-S:aHqg}).  Hence, any of the last three
  conditions give that the manifold is hyperHermitian and $\partial_I
  (\omega_J + i\omega_K) = 0$ and we obtain HKT from Grantcharov \& Poon.

  Alternatively, we may prove the result just using tools contained in the
  present paper.  Suppose $N_J = 0$ and $Jd\omega_J=Kd\omega_K$.  Then
  $Jd\omega_J$ and hence $Kd\omega_K$ lie in $\Wc{3+4,J}$.  However,
  Proposition~\ref{prop:naddd}\itref{item:naddd-4} gives that $Kd\omega_K +
  Id\omega_I \in \Wc{3+4,J}$, so we have $Id\omega_I \in \Wc{3+4,J}$ too.
  Now let us use Proposition~\ref{prop:naddd}\itref{item:naddd-2} for the
  integrability of~$J$.  We have
  $0=K_{(23)}(-J_{(12)}-J_{(13)}+J_{(23)}-1)(Kd\omega_K-Id\omega_I)$.  But
  $J_{(12)}+J_{(23)}+J_{(13)}=1$ on~$\Wc{3+4,J}$ and $Jd\omega_J =
  Kd\omega_K$, so
  \begin{equation*}
    J_{(23)} (Jd \omega_J - Id\omega_I) = Jd\omega_J - Id\omega_I.
  \end{equation*}
  Skew-symmetrising both sides of the identity, we find that $Jd\omega_J -
  Id\omega_I = 3 (Jd\omega_J - Id\omega_I)$.  So, $Id\omega_I = Jd\omega_J
  = Kd\omega_K$.
\end{proof}

Next we describe the very special situation for four-dimensional
HKT-manifolds.

\begin{proposition}
  If $M$ is an almost hyperHermitian $4$-manifold, then the following
  conditions are equivalent:
  \begin{enumerate}
  \item $M$ is an HKT-manifold;
  \item the three Lee one-forms are equal, i.e., $Id^*\omega_I =
    Jd^*\omega_J = Kd^*\omega_K$;
  \item\label{item:HKT-3} the almost complex structures $I$ and $J$ are
    integrable;
  \item\label{item:HKT-4} the almost Hermitian structures corresponding to
    $I$ and $J$ are locally conformally K\"ahler, so $M$ is locally
    conformally hyperK\"ahler.
  \end{enumerate}
\end{proposition}

\begin{proof}
  For dimension $4$, the Gray-Hervella modules $\Wc1$ and $\Wc3$ are zero,
  we have $I \theta \wedge \omega_I = J \theta \wedge \omega_J = K \theta
  \wedge \omega_K$, for all one-forms~$\theta$, and any three-form may be
  written in this way.  If $M^4$ is an HKT-manifold, we see that the almost
  Hermitian structures are of type $\Wc4$ and that
  \begin{equation*}
    T = Ad\omega_A = - At \wedge \omega_A,
  \end{equation*}
  where $t = A\Lambda_A T = - \Lambda_A d\omega_A = A d^* \omega_A$.  On
  the other hand, if the three Lee forms are equal to a one-form~$t$, then
  \begin{equation*}
    Id\omega_I = - It \wedge \omega_I = - Jt \wedge \omega_J = Jd\omega_J.
  \end{equation*}
  Hence $Id\omega_I = Jd\omega_J = Kd\omega_K$ and $M$ is HKT.

  For conditions \itref{item:HKT-3} or \itref{item:HKT-4}, the three almost
  complex structures are integrable, so the almost Hermitian structures
  have a common Lee form, by~\cite{Cabrera-S:aHqg}.
\end{proof}

In \S\ref{sec:exterior} it was shown that for any almost
quaternion-Hermitian manifold, the exterior derivatives of the three local
K\"ahler forms of an adapted basis $I$, $J$, $K$ satisfy the
identities~\eqref{eq:ijk}.  When the manifold is HKT, additional identities
are also satisfied.

\begin{lemma}
  \label{lem:hktonef}
  For a $4n$-dimensional HKT-manifold, the exterior derivatives
  $d\omega_I$, $d\omega_J$ and $d\omega_K$ satisfy
  \begin{equation}
    \label{eq:HKT-d}
    t = - \Lambda_I d\omega_I = K\Lambda_J d\omega_I = - J\Lambda_K
    d\omega_I, \detc, 
  \end{equation}
  where $t = Id^*\omega_I = Jd^*\omega_J = Kd^*\omega_K$.  Furthermore,
  $I\lambda_I = J\lambda_J = K\lambda_K = \tfrac1{2n} t$, $\theta^\xi =
  \theta_I^\xi = \theta_J^\xi = \theta_K^\xi$ and the one-forms
  corresponding to the $E$-parts of~$\beta_A$ are such that
  \begin{gather}
    \label{eq:nu-4-t}
    \nu^I_4 = \nu^J_4 = \nu^K_4 = 2t,\\
    \label{eq:nu-3-t}
    \nu^I_3 =\nu^J_3 = \nu^K_3 = 32 \theta^\xi = \tfrac4{2n+1} t,
  \end{gather}
  where the second line holds for $n>1$.
\end{lemma}

\begin{proof}
  Since $2T = \beta_I = \beta_J = \beta_K$, we have
  \begin{equation*}
    I\Lambda_I \beta_A = J\Lambda_J\beta_A = K\Lambda_K\beta_A,
  \end{equation*}
  from which we obtain \eqref{eq:HKT-d} and~\eqref{eq:nu-4-t},
  via~\eqref{eq:nu-I4}.  Now, using equation~\eqref{eq:beta-theta}, we have
  $(2n+1) (n-1) \nu^A_3 = 4 (n-1) A d^* \omega_A$ and
  hence~\eqref{eq:nu-3-t}.
\end{proof}

\section{Quaternion-K\"ahler manifolds with torsion}
\label{sec:QKT}

A genuinely quaternionic analogue of HKT geometry also arises in the
physics literature via the theory of super-symmetric sigma models.  In this
section we give a definition in terms of intrinsic torsion, relate this
definition to the existence of connections with skew-symmetric torsion,
provide different characterisations of the geometry and describe the
relationship with HKT geometry.  Important mathematical work in this
direction was previously done by Ivanov~\cite{Ivanov:QKT}.  Here we
concentrate on the intrinsic geometry, fit the geometry into our general
formalism and improve and clarify a number of his results.

\begin{definition}
  An almost quaternion\bdash Hermitian manifold of dimension $4n\geqslant8$
  is \emph{QKT} (\emph{quaternion-K\"ahler with torsion}) if its intrinsic
  torsion lies in $(K+E)H$.
\end{definition}

As in other cases, we may write this condition on the intrinsic torsion in
terms of three-forms.

\begin{lemma}
  The intrinsic torsion~$\xi$ lies in $(K+E)H$ precisely when it is given
  by a three-form~$\psi\in (K+E)H\subset\Lambda^3T^*M$ via
  \begin{equation}
    \label{eq:char-KH-EH}
    \inp Y{\xi_X Z} = \Bigl(3\psi + \sum_{A=I,J,K} \bigl(-A_{(23)}\psi + 
    \tfrac2n A \theta^\psi \otimes \omega_A\bigr) \Bigr) (X,Y,Z), 
  \end{equation}
  where $\theta^\psi = I\Lambda_I\psi = J\Lambda_J\psi = K\Lambda_K\psi$.
  Moreover, for a given intrinsic torsion $\xi\in(K+E)H$ we have that
  $\psi\in\Lambda^3T^*M$~is unique and given by
  \begin{equation*}
    48(n-1)\psi  = (n-1)d^*\Omega + \tfrac12 \sum_{A= I,J,K} A\theta^{d^*
    \Omega} \wedge\omega_A, 
  \end{equation*}
  where $\Omega$ is the fundamental four-form~\eqref{eq:four-form} and
  $d^*$~is the co-derivative.
\end{lemma}

\noindent
When applying this result it is often useful to recall the
formula~\cite{Cabrera:aqh}
\begin{equation*}
  d^*\Omega = 2 \sum_{A=I,J,K} (d^* \omega_A \wedge \omega_A - A d
  \omega_A). 
\end{equation*}

Let us now demonstrate how the QKT condition relates to connections with
skew-symmetric torsion and so the original definition of Howe, Opfermann \&
Papadopoulos \cite{Howe-OP:QKT}.  Recall that $(K+E)H\subset\Lambda^3T^*M$
is the $(+3)$-eigenspace of the operator~$\mathcal L$ given
in~\eqref{eq:L}.

\begin{theorem}
  \label{thm:QKT}
  An almost quaternion\bdash Hermitian manifold $M$ is QKT if and only if
  there exists a metric connection $\Nqkt = \LC + \frac12 T$ that is
  quaternionic and whose $(0,3)$-torsion $T(X,Y,Z) = \inp X{T(Y,Z)}$ is a
  three-form satisfying in $\mathcal L T=-3T$.  When $M$~is QKT, $\Nqkt$~is
  the unique $\SP(n)\SP(1)$-connection on $M$ with skew-symmetric torsion.
\end{theorem}

Concretely, we claim that the intrinsic torsion~$\xi$ of the QKT structure
is given by~\eqref{eq:char-KH-EH} with $\psi=-\frac18T$ and that so
\begin{equation*}
  T = - \tfrac16 d^* \Omega - \tfrac1{12(n-1)} \sum_{A=I,J,K}
  A \theta^{d^* \Omega} \wedge \omega_A.
\end{equation*}

Note that in \S\ref{sec:examples} we will provide examples of almost
quaternion\bdash Hermitian manifolds that are not QKT but none-the-less
admit $\SP(n)\SP(1)$\bdash connections with skew-symmetric torsion.  

\begin{proof}
  If $\xi \in (K+E)H$, then $\xi$ is given by
  equation~\eqref{eq:char-KH-EH} for some~$\psi$ in~$(K+E)H$.  Putting $T =
  -8\psi$ we find that $\Nqkt = \LC + \frac12 T$ is metric and,
  via~\eqref{eq:aqH-torsion}, quaternionic.

  Conversely, if $M$ has such a connection~$\Nqkt$, then
  \begin{equation}
    \label{eq:QKT-omega}
    \LC \omega_I = \gamma_K \otimes \omega_J - \gamma_J \otimes \omega_K -
    \tfrac12 (I_{(2)} + I_{(3)}) T,
  \end{equation}
  where $\gamma_A$, $A=I,J,K$, are the one-forms given
  by~\eqref{eq:q-nabla} for $\Nt=\Nqkt$.  Using equation~\eqref{eq:LC-xi},
  we find
  \begin{equation}
    \label{eq:T-pre}
    \begin{split}
      \inp Y{\xi_XZ} &= -\frac18 \Bigl(3T-\sum_{A=I,J,K}
      A_{(23)}T\Bigr)(X,Y,Z)\eqbreak +\frac12\sum_{A=I,J,K}
      (\lambda_A-\gamma_A)\otimes\omega_A(X,Y,Z).
    \end{split}
  \end{equation}
  As $\xi \in T^*M \otimes \Lambda_0^2 E \sym2 H$, we have
  $\inp{Ae_i}{\xi_Xe_i} =0$ and find
  \begin{equation}
    \label{eq:lambda-gamma}
    \lambda_A - \gamma_A = - \tfrac1{2n} A t,
  \end{equation}
  with $t=I\Lambda_IT=J\Lambda_JT=K\Lambda_KT$.  Using
  \eqref{eq:lambda-gamma} equation~\eqref{eq:T-pre}, we obtain
  equation~\eqref{eq:char-KH-EH} with $\psi = - \tfrac18 T \in (K+E)H$.
\end{proof}

\begin{remark}
  The situation for $4$-dimensional almost quaternion-Hermitian manifolds
  is very special.  Here the Levi-Civita connection is always quaternionic,
  i.e.,
  \begin{equation*}
    \LC I = \lambda_K \otimes J - \lambda_J \otimes K, \detc
  \end{equation*}
  Also in this dimension, we have $\Lambda^3T^*M\cong T^*M$ and any
  three-form $T$ may be written as $T=-At\wedge\omega_A$ for some
  one-form~$t$ valid for $A=I$, $J$ and~$K$.  In this way, given any
  $t\in\Omega^1(M)$, we may construct $\Nt=\LC+\tfrac12 T$ and find
  \begin{equation*}
    \Nt I = \gamma_K \otimes J - \gamma_J \otimes K, \detc,
  \end{equation*}
  where $\gamma_A = \lambda_A + \tfrac12 At$.  Hence $\Nt$ is a connection
  with skew-symmetric torsion preserving the almost quaternion-Hermitian
  structure.  However, in this case, $\Nt$ is not unique.
\end{remark}

Forgetting the metric of an almost quaternion-Hermitian structure we are
left with an almost quaternionic structure.  This is an \emph{integrable}
quaternionic structure if there is a torsion-free quaternionic
connection~$\Nq$, i.e., $\Nt=\Nq$ satisfies~\eqref{eq:q-nabla}; this is a
weaker condition than integrability of $I$, $J$ and~$K$.  In the presence
of a compatible metric, integrability of the quaternionic structure is
equivalent to the vanishing of $\xi_{S^3H}$, the
$(\Lambda^3_0E+K+H)S^3H$-part of the intrinsic torsion,
cf.~\cite{Salamon:quaternionic}.  In dimension four, this condition is just
self-duality of the conformal structure.  In dimension~$8$, the module
$\Lambda^3_0E$ is zero, so $\xi\in (K+E)(S^3H+H)$ and integrability implies
that $\xi\in(K+E)H$.  We thus have:

\begin{proposition}
  Any compatible metric on an eight-dimensional quaternionic manifold is
  QKT.\qed
\end{proposition}

This applies for example to any metric compatible with Joyce's invariant
hypercomplex structure on~$\SU(3)$ \cite{Joyce:hypercomplex}.

The one-form $t$ in the proof of Theorem~\ref{thm:QKT} has independent
importance.

\begin{definition}
  For a QKT manifold with torsion three-form~$T$ the \emph{torsion
  one-form} $t$ is defined by
  \begin{equation*}
    t = I\Lambda_IT,
  \end{equation*}
  for any compatible almost complex structure~$I$.
\end{definition}

There are many alternative expressions for~$t$:

\begin{lemma}
  \label{lem:one-form}
  For a $4n$-dimensional $QKT$-manifold, $n>1$, the torsion one-form $t$
  satisfies
  \begin{equation*}
    - \tfrac{3(n-1)}{4n} t
    = \left(\xi_{e_i}e_i\right)^\flat
    = -\tfrac32 A\eta_A
    = \tfrac1{16n} \ast (\ast d\Omega \wedge \Omega)
    = -\tfrac3{8n} A\Lambda_A d^*\Omega ,
  \end{equation*}
  where $\eta_A$ is given by~\eqref{eq:eta-I} and $X^\flat=\inp X\cdot$.
\end{lemma}

\begin{remark}
  \label{rem:EH}
  The one-form
  \begin{equation*}
    \ast ( \ast d \Omega \wedge \Omega) = -2 \sum_{A=I,J,K} A\Lambda_A d^* \Omega
  \end{equation*}
  was considered in~\cite{Cabrera:aqh} in relation the $EH$-component
  of~$\xi$ in general.  One finds that $ \ast ( \ast d \Omega \wedge
  \Omega) = 16 n \left( \xi_{e_i} e_i \right)^\flat$.  There it was noted
  that
  \begin{equation*}
    4n ( J \eta_J + K \eta_K) = 2 I\Lambda_I d^* \Omega = - \ast \left(
      \ast d \Omega \wedge \omega_I \wedge \omega_I \right).
  \end{equation*}
\end{remark}

\begin{proof}
  Using~\eqref{eq:xi-alpha}, one finds $\left(\xi_{e_i} e_i \right)^\flat =
  - \tfrac12 \sum_{A=I,J,K} A \eta_A$.  When $\xi$ is in~$(K+E)H$, we have
  $I \eta_I = J \eta_J = K \eta_K$.  Therefore, $\left( \xi_{e_i} e_i
  \right)^\flat = - \frac32 I \eta_I$.  On the other hand, using
  $\Psi=-\frac18 T$ in equation~\eqref{eq:char-KH-EH}, we obtain $\left(
    \xi_{e_i} e_i \right)^\flat = - \frac{3(n-1)}{4n} t$.

  The remaining equalities follow from Remark~\ref{rem:EH}
\end{proof}

Further relations between the torsion three-form $T$ and the four-form
$\Omega$ follow from equation~\eqref{eq:QKT-omega}.  In particular, we have
\begin{equation*}
    \LC \Omega = - \sum_{A=I,J,K} (A_{(2)}+A_{(3)}) T \wedge
    \omega_A, \quad d \Omega = - 2 \sum_{A=I,J,K} AT \wedge
    \omega_A.
\end{equation*}

Let us now give a characterisation of QKT manifolds.

\begin{theorem}
  \label{thm:QKT1}
  An almost quaternion-Hermitian manifold~$M$ is QKT if and only if for
  each local adapted basis $I$, $J$, $K$ there are local one-forms
  $\gamma_I$, $\gamma_J$, $\gamma_K$ such that
  \begin{equation}
    \label{eq:omega-diff}
    \begin{split}
      \beta_J - \beta_I &= I d \omega_I - J d \omega_J\\
      &= I (K \gamma_K) \wedge \omega_I - J(K \gamma_K) \wedge \omega_J + K
      \gamma_K^- \wedge \omega_K, \detc,
    \end{split}
  \end{equation}
  where $\gamma_K^- = I \gamma_I - J \gamma_J$;

  In this case, the skew-symmetric torsion three-form $T$ is given by
  \begin{align}
    \label{eq:T-Id}
    2T
    &= 2(I d \omega_I + J (K \gamma_K) \wedge \omega_J + K (J
    \gamma_J) \wedge \omega_K)\\
    \label{eq:T-beta}
    & = \beta_I + J (I \gamma_I) \wedge \omega_J + K (I \gamma_I) \wedge
    \omega_K + I \gamma_I^+ \wedge \omega_I, \detc
  \end{align}
  where $\gamma_I^+ = J \gamma_J + K \gamma_K$ and 
  \begin{equation}
    \label{eq:gammaI}
    2 (n-1) I \gamma_I  =  (\Lambda_J + I\Lambda_K  ) d \omega_J 
    = (\Lambda_K - I\Lambda_J)d\omega_K, \detc.
  \end{equation}
\end{theorem}

\begin{remark}
  Equation~\eqref{eq:T-beta} also implies that, for dimensions strictly
  greater than four, the QKT-connection~$\Nqkt$ is unique.  This fact was
  already proved by Ivanov, who also characterised QKT-manifolds by the
  differences $\left( I d \omega_I \right)_{\Wc{3+4,I}} - \left( J d
    \omega_J \right)_{\Wc{3+4,J}}$ \cite[Theorem 2.2]{Ivanov:QKT}.  Our
  Theorem~\ref{thm:QKT1} can be considered as an improved version, based on
  the three-forms $\beta_A$, which are automatically in~$\Wc{3+4,A}$.
\end{remark}

\begin{remark}
  Let us write $\Wm^{(E)}_{3I} (\eta_I)$ and $\Wm_{4I} (\eta_I)$ for the
  right-hand sides of equations \eqref{eq:beta-E3} and~\eqref{eq:beta-E4},
  respectively.  Equation~\eqref{eq:T-beta} can then be written
  \begin{equation*}
    2 T = \beta_A - \Wm^{(E)}_{3A} \left( 2 A \gamma A\right) - \Wm_{4A}
    \left( 2 A\gamma_A + (2n-1) \gamma_A^+ \right).
  \end{equation*}
  Consequently, we see that QKT-manifolds have
  \begin{gather*}
    \beta_A^{(3)}=0, \quad
    \beta_I^{(K)} = \beta_J^{(K)} = \beta_K^{(K)},\\
    \nu^A_3 = 2 A \gamma_A + \tfrac4{2n+1} t \quad\text{and}\quad \nu^A_4 =
    2 A \gamma_A + (2n-1) \gamma_A^+ + 2 t,
  \end{gather*}
  using equation~\eqref{eq:lambda-gamma} and Lemma \ref{lem:one-form}.
  This should be compared with the results of Lemma~\ref{lem:hktonef}.
\end{remark}

\begin{proof}
  Suppose \( M \) is QKT.  From $\Nqkt = \LC + \frac12 T$ and
  equation~\eqref{eq:QKT-omega}, we have
  \begin{equation}
    \label{eq:d-omegaI}
    d\omega_I = - IT + \gamma_K\wedge\omega_J - \gamma_J\wedge\omega_K.
  \end{equation}
  Multiplying by~$I$ gives \eqref{eq:T-Id} and
  equation~\eqref{eq:omega-diff} follows.

  Conversely, if equation~\eqref{eq:omega-diff} is satisfied for some
  one-forms $\gamma_A$, then we may consider the three-form~$T$ given
  by~\eqref{eq:T-Id} and use \eqref{eq:omega-diff} to obtain the
  alternative expression~\eqref{eq:T-beta}, see that it is unchanged by a
  cyclic permutation of $I,J,K$, and hence globally defined.  By
  Proposition~\ref{prop:naddd}\itref{item:naddd-4}, we find that $T \in
  \Wc{3+4,A}$, for $A=I,J,K$.  Hence $T \in (K+E)H$ and $M$~is a
  $QKT$-manifold for the connection $\Nqkt = \LC + \frac12 T$.

  Finally, using equation~\eqref{eq:d-omegaI}, we get
  \begin{equation*}
    K\Lambda_J d \omega_I = t + J \gamma_J
    + (2n-1) K\gamma_K\quad\text{and}\quad
    - \Lambda_I d \omega_I = t + J\gamma_J + K \gamma_K.
  \end{equation*}
  Using the identity~\eqref{eq:ijk}, now provides the claimed expressions
  for~$I\gamma_I$.
\end{proof}

The contraction identity~\eqref{eq:ijk} used above is valid for all almost
quaternion\bdash Hermitian manifolds.  When the structure is QKT, there are
additional identities of this type.

\begin{proposition}
  \label{prop:QKT2}
  For a local adapted basis $I,J,K$ of a QKT-manifold, the following
  identities are satisfied:
  \begin{gather}
    \label{eq:QKT-coder}
    Id^*\omega_I = t + J\gamma_J + K\gamma_K = J\lambda_J + K\lambda_K +
    2I\eta_I, \detc, \\ 
    J\gamma_J - K\gamma_K = J\lambda_J - K\lambda_K = -(Jd^*\omega_J -
    Kd^*\omega_K),\detc, \label{eq:gamma-diff} \\
    J\Lambda_Kd\omega_I + K\Lambda_Id\omega_I = -2(n-1) (\Lambda_Jd\omega_J
    - \Lambda_Kd\omega_K), 
    \detc,  \label{eq:newddd1} \\
    J\Lambda_Id\omega_K + K\Lambda_Id\omega_J = (2n-1) (\Lambda_Jd\omega_J
    - \Lambda_Kd\omega_K), 
    \detc \label{eq:newddd2}
  \end{gather}
\end{proposition}

\begin{proof}
  As $M$ is a QKT-manifold, we have from equation~\eqref{eq:QKT-omega}
  \begin{equation*}
    d^* \omega_I = - (\nabla_{e_i}\omega_I)(e_i,\cdot) = -I t - I(J\gamma_J
    + K\gamma_K)   
  \end{equation*}
  giving the first equality of equation~\eqref{eq:QKT-coder}.  On the other
  hand, it was shown in \cite{Cabrera-S:aHqg} that $I d^* \omega_I = J
  \lambda_J + K \lambda_K + J \eta_J + K \eta_K$.  But, for manifolds of
  type $KH+EH$, we have $I \eta_I = J \eta_J = K \eta_K$ by Corollary
  \ref{cor:E}, and we obtain the second equality of
  equation~\eqref{eq:QKT-coder}.

  Equation~\eqref{eq:gamma-diff} is an immediate consequence of
  equation~\eqref{eq:QKT-coder}.

  Finally, equation~\eqref{eq:newddd1} and equation~\eqref{eq:newddd2} can
  be deduced from equations~\eqref{eq:gammaI}, \eqref{eq:ijk}
  and~\eqref{eq:gamma-diff}, using $Id^* \omega_I = - \Lambda_I d\omega_I$.
\end{proof}

\begin{remark}
  It is not hard to prove that an almost quaternion-Hermitian
  $4n$-manifold, $n>1$, is of type $\Lambda^3_0 E \sym3 H + K \sym3 H +
  \Lambda^3_0 EH + KH+ EH$ if and only if equation~\eqref{eq:newddd1} is
  satisfied for any local adapted basis $I$, $J$, $K$.  In other words,
  \eqref{eq:newddd1} characterises the vanishing of the $E\sym3H$-component
  of the intrinsic torsion.
\end{remark}

Let us now turn to the question of integrability of compatible almost
complex structures $I$, $J$ and~$K$.  Taking the following identity
\begin{equation*}
  N_I(X,Y) = - (\LC_X I)IY - (\LC_{IX}I)Y + (\LC_YI)IX + (\LC_{IY}I)X,
\end{equation*}
into account, we obtain that the corresponding $(0,3)$-tensor $N_I =
\inp\cdot{N_I(\cdot,\cdot)}$ is given by
\begin{equation}
  \label{eq:NI}
  N_I = J\gamma_I^- \wedge \omega_J - K\gamma_I^- \wedge \omega_K
  - J\gamma_I^- \otimes \omega_J + K\gamma_I^- \otimes \omega_K, \detc
\end{equation}
Note that equation~\eqref{eq:lambda-gamma} gives $\gamma_A^- = \lambda_A^-
$, where $\lambda_I^-= J \lambda_J - K \lambda_K$.

Fixing the almost complex structure $I$ and under the action of the
subgroup $\SP(n) \Un(1)$ of $\Un(2n)_I$, it was noted in
\cite{Cabrera-S:aHqg} that, for $n>1$, $\Wc1 \otimes \mathbb C =
(\Lambda_0^3 E + E) (L^3 + \bar L^3)$ and $\Wc2 \otimes \mathbb C = (K + E)
(L^3 + \bar L^3)$.  As in \S~\ref{sec:beta}, we get
\begin{equation*}
  \Wc1 = [\Lambda_0^3 EL^3]_{\mathbb R} + [EL^3]_{\mathbb R1}, \qquad
  \Wc2 = [KL^3]_{\mathbb R} + [EL^3]_{\mathbb R2}.
\end{equation*}
It is well known that, in general, the tensor~$N_I$ belongs to~$\Wc{1+2}$.
In our case, equation~\eqref{eq:NI} gives us that $N_I \in [EL^3]_{\mathbb
R1} + [EL^3]_{\mathbb R2}$, i.e., the components of~$N_I$ in $[\Lambda_0^3
EL^3]_{\mathbb R}$ and $[KL^3]_{\mathbb R}$ vanish.

Using equations \eqref{eq:NI} and~\eqref{eq:gamma-diff}, we thus have:

\begin{corollary}
  \label{cor:QKT3}
  Let $M$ be a QKT-manifold of dimension $4n>4$.  Then, for any local
  adapted basis $I,J,K$, the following conditions are equivalent:
  \begin{enumerate}
  \item the almost complex structure $I$ is integrable;
  \item the Lee forms $J d^* \omega_J$ and $K d^* \omega_K$ are equal;
  \item the one-forms $J \lambda_J$ and $K \lambda_K$ are equal;
  \item the one-forms $J \gamma_J$ and $K \gamma_K$ are equal;
  \item the equation $ \Lambda_K d \omega_J = - \Lambda_J d \omega_K$ is
    satisfied; 
  \item the equation $ J\Lambda_K d\omega_I = - K\Lambda_J d \omega_I$ is
    satisfied;
  \item the equation $ J\Lambda_I d \omega_K = - K\Lambda_I d \omega_J$ is
    satisfied. 
    \qed
  \end{enumerate}
\end{corollary}

\begin{corollary}
  Let $M$ be a $4n$-dimensional, $(n>1)$, almost quaternion-Hermitian with
  a global adapted basis $I,J,K$, i.e., $M$ is equipped with an almost
  hyperHermitian structure.  Then the following conditions are equivalent
  \begin{enumerate}
  \item $M$ is an HKT-manifold;
  \item $M$ is a $QKT$-manifold such that $I \lambda_I = J \lambda_J = K
    \lambda_K = \frac1{2n} t$;
  \item $M$ is a $QKT$-manifold such that $I d^* \omega_I = J d^* \omega_J
    = K d^* \omega_K$ and
    \begin{equation*}
      I\Lambda_K d\omega_J  =  Id^* \omega_I, \detc
    \end{equation*}
  \end{enumerate}
\end{corollary}

\begin{proof}
  This is an immediate consequence of Theorem~\ref{thm:QKT1},
  Proposition~\ref{prop:QKT2} and Corollary~\ref{cor:QKT3}.
\end{proof}

Finally, lets us look at the two special types of QKT-manifolds with
intrinsic torsion in one summand of~$(K+E)H$.

\begin{lemma}
  \label{lem:aqh-KH}
  An almost quaternion-Hermitian manifold $M$~of dimension $4n>4$ is of
  type $\Lambda_0^3 E \sym3 H + K \sym3 H + \Lambda_0^3 E H + K H$ if and
  only if, for any local adapted basis $I$, $J$, $K$, we have
  \begin{equation}
    \label{eq:aqh-KH}
    -I\Lambda_K d \omega_J = (n-1) I d^* \omega_I - n J d^* \omega_J -
    (n-1) K d^* \omega_K, \detc 
  \end{equation}
\end{lemma}

\begin{proof}
  We have $\theta^\xi = \theta^\xi_I = \theta^\xi_J = \theta^\xi_K =0$.
  Equations \eqref{eq:beta-3A} and~\eqref{eq:beta-4A} then give
  \begin{equation*}
    \nu^A_3 = 2 A \lambda_A, \qquad \nu^A_4 = (2n-1)
    \lambda_A^+ + 2 A \lambda_A.
  \end{equation*}
  From these equalities together with equations \eqref{eq:beta-E4},
  \eqref{eq:beta-theta} and~\eqref{eq:ladd2}, we deduce
  equation~\eqref{eq:aqh-KH}.
\end{proof}

\begin{theorem}
  Let $M$ be an almost quaternion-Hermitian $4n$-manifold, $n>1$.
  \begin{enumerate}
  \item $M$ is of type $KH$ if and only if $M$ is a $QKT$-manifold and, for
    any local adapted basis $I$, $J$, $K$, equation~\eqref{eq:aqh-KH} holds.
  \item $M$ is of type $EH$ if and only if there exists a
    global one-form $t$ defined on $M$ such that, for any local adapted
    basis $I$, $J$, $K$, we have
    \begin{equation}
      \label{eq:aqh-EH}
      \begin{split}
        d \omega_I & =  - \tfrac1{2n+1} t \wedge \omega_I - K \left( K
          \lambda_K + \tfrac1{2n(2n+1)} t \right)  \wedge \omega_J
        \eqbreak + J \left( J \lambda_J + \tfrac1{2n(2n+1)} t \right) \wedge
        \omega_K, \detc \nonumber
      \end{split}
    \end{equation}
    In this case, $M$ is called a locally conformal quaternionic K\"ahler
    manifold and the one-form $t$ is given by
    \begin{equation*}
      2(n-1) t = 2n I d^* \omega_I - K\Lambda_J d \omega_I + J\Lambda_K d
      \omega_I, \detc 
    \end{equation*}
  \end{enumerate}
\end{theorem}

\begin{proof}
  The first part follows directly from Lemma~\ref{lem:aqh-KH}.

  For $M$ of type~$EH$, then $ T = - \frac1{2n+1} \sum_{A=I,J,K} A t \wedge
  \omega_A $.  Equation~\eqref{eq:aqh-EH} then follows from equations
  \eqref{eq:d-omegaI} and~\eqref{eq:QKT-coder} and
  Lemma~\ref{lem:one-form}.
\end{proof}

\section{Almost Hermitian structures}
\label{sec:Gray-Hervella} 

In this section, we will consider almost Hermitian manifolds $M$ of
dimension $2n$ and the classification of Gray \& Hervella \cite{Gray-H:16}.
These manifolds are equipped with an almost complex structure~$I$
compatible with a Riemannian metric $\inp\cdot\cdot$.  Therefore, their
orthogonal frame bundles can be reduced to the unitary group~$\Un(n)$.

By identifying the intrinsic $\Un(n)$-torsion $\IaH$ with
$\LC\omega_I$ via $\IaH \mapsto - \IaH \omega_I = \LC \omega_I
$, Gray \& Hervella \cite{Gray-H:16} gave conditions characterising classes
of almost Hermitian manifolds by means of the covariant
derivative~$\LC\omega_I$.  The space of intrinsic $\Un(n)$-torsions is then
isomorphic to the space~$\Wc{}\cong T^*M\otimes\Lambda^{\{2,0\}}$ of
covariant derivatives of the K\"ahler form~$\omega_I$.  Under the action of
$\Un(n)$, $n >2$, $\Wc{}$~decomposes into four irreducible modules,
\begin{equation*}
  \Wc{} = \Wc1 + \Wc2 + \Wc3 + \Wc4 \cong
  \Lambda^{\{3,0\}} + [U^{3,0}]_{\mathbb R} + \Lambda^{\{2,1\}}_0
  + \Lambda^{\{1,0\}}.
\end{equation*}
Therefore, for $n>2$, there are $2^4=16$ classes of almost Hermitian
manifolds.  For $n=2$, $\Wc1 = \Wc3 =\{0\}$ and there are only $4$ classes.

On the other hand, because $d \omega_I \in \Lambda^3 T^*M = \Wc{1+3+4}$,
only partial information about $\IaH$ can be recovered from the
exterior derivative~$d \omega_I$.  The remaining component can be found in
the Nijenhuis $(0,3)$-tensor $N_I \in \Wc{1+2} \subset T^*M \otimes
\Lambda^2 T^*M$, which is often more convenient to work with than~$\LC$.
Table~\ref{tab:ahtypes} lists conditions characterising the classes of
almost Hermitian manifolds in terms of tensors $d\omega_I$ and~$N_I$.  The
symbol $\mathcal N_I$ denotes the three-form obtained by
skew-symmetrisation of $N_I$, i.e., $\mathcal N_I (X,Y,Z) = \sumcic_{XYZ}
N_I (X,Y,Z)$.  The conditions are found by studying the $\Un(n)$-maps
\begin{gather*}
  \IaH \mapsto -  \sumcic_{XYZ} (\IaH_X \omega_I)(Y,Z)
  = d \omega_I (X,Y,Z),\\
  \begin{split}
    \IaH &\mapsto - ( \IaH_Z \omega_I ) (IX,Y) - (
    \IaH_{IZ} \omega_I ) (X,Y)\eqbreak - ( \IaH_{Y} \omega_I )
    (IZ,X) - ( \IaH_{IY} \omega_I ) (Z,X)
  \end{split}
\end{gather*}
and recalling the following well-known identity \cite{Gray:minimal}
\begin{equation*}
  \begin{split}
    N_I (X,Y,Z) &= (\LC_Z \omega_I ) (IX,Y) + ( \LC_{IZ} \omega_I ) (X,Y)
    \eqbreak + ( \LC_{Y} \omega_I ) (IZ,X) + ( \LC_{IY} \omega_I ) (Z,X).
  \end{split}
\end{equation*}

\begin{table}[tp]
  \centering
  \begin{tabular}{ll}
    \toprule
    $\mathcal K$ &
    $N_I=0$ and $d \omega_I=0$ \\ 
    $\Wc1 = \mathcal {NK}$ &
    $d \omega_I = - \frac34 IN_I$ \\
    $\Wc2 =\mathcal{AK} $ &
    $d \omega_I=0$ \\
    $\Wc3$ &
    $N_I=0$ and $d^* \omega_I =0 $ \\
    $\Wc4 = \mathcal{LCK}$ &
    $N_I=0$ and $d \omega_I = - \frac1{n-1} I d^* \omega_I \wedge \omega_I $ \\ 
    $\Wc{1+2}$ &
    $d \omega_I = - \frac14 I \mathcal N_I$ \\
    $\Wc{1+3}$ &
    $\mathcal N_I=3N_I$ and $d^* \omega_I=0$ \\
    $\Wc{1+4}$ &
    $d \omega_I = - \frac34 I N_I - \frac1{n-1} I d^* \omega_I \wedge
    \omega_I$ \\ 
    $\Wc{2+3}$ &
    $\mathcal N_I=0$
    and $d^* \omega_I=0$ \\
    $\Wc{2+4}$ &
    $d \omega_I = - \frac1{n-1} I d^* \omega_I \wedge \omega_I $ \\
    $\Wc{3+4}$ &
    $N_I =0$ \\
    $\Wc{1+2+3}$ &
    $ d^* \omega_I =0 $ \\
    $\Wc{1+2+4}$ & 
    $d \omega_I = - \frac14 I \mathcal N_I - \frac1{n-1} I d^* \omega_I \wedge \omega_I$ \\
    $\Wc{1+3+4}$ & 
    $\mathcal N_I = 3 N_I$ \\
    $\Wc{2+3+4}$ &
    $\mathcal N_I = 0$ \\
    $\Wm$ &
    no relation \\
    \bottomrule
  \end{tabular}
  \caption{The Gray-Hervella classes of almost Hermitian structures
  characterised by the Nijenhuis tensor and the exterior derivative of the
  K\"ahler form.}
  \label{tab:ahtypes}
\end{table}

\section{Twisting}
\label{sec:twist}

In this section we consider the effects of \enquote{twisting} an almost
hyperHermitian in the sense of~\cite{Swann:T} where this construction was
shown to be an interpretation of T-duality.

Let $(M,I,J,K,g=\inp\cdot\cdot)$ be an almost hyperHermitian manifold.
Suppose that $X$ is a tri-holomorphic isometry, so
\begin{equation*}
  L_Xg = 0, \qquad L_XI=0=L_XJ=L_XK.
\end{equation*}
Let $F_\theta$ be a closed $2$-form with $L_XF_\theta=0$ and choose a
nowhere vanishing function $a\in C^\infty(M)$ so that $X^\theta = X\hook
F_\theta=-da$.  If $P\to M$ is a principal $\mathbb R$-bundle with
connection $\theta$ whose curvature is $F_\theta$, then $X$ may be lifted
to a transverse vector field $\tilde X$ on~$P$ that preserves $\theta$: the
vertical component is given by $aY$ where $Y$ generates the principal
action.  Locally the \emph{twist}~$W$ of $M$ by $(X,F_\theta,a)$ is the
quotient $W=P/\langle \tilde X\rangle$ with the geometry induced from the
horizontal space~$\mathcal H=\ker\theta$.  Each $X$-invariant
$(0,p)$-tensor~$\kappa$ on~$M$ is \emph{$\mathcal H$-related} to a unique
$(0,p)$-tensor $\kappa^W$ on~$W$ defined by the condition that the
pull-backs to~$P$ agree on~$\mathcal H$.  Exterior differentiation on~$W$
then corresponds to the twisted derivative $d^W =
d-F_\theta\wedge\tfrac1aX\hook$ on invariant forms on~$M$.

If we twist an almost hyperHermitian structure the three-form
$\beta_I=Jd\omega_J+Kd\omega_K$ transforms to
\begin{equation*}
  \beta_I^W = \beta_I - \tfrac1aX^\flat\wedge (J+K)F_\theta.
\end{equation*}
Thus to compute the change of the intrinsic torsion we should decompose
$\beta_I^W$ as in~\eqref{eq:beta}.  We start by decomposing $F_\theta \in
\Lambda^2T^*M = S^2H + \Lambda^2_0ES^2H + S^2E$:
\begin{equation*}
  F_\theta = (\mu_I\omega_I + \mu_J\omega_J + \mu_K\omega_K) +
  (I_{(1)}\kappa_I + J_{(1)}\kappa_J + K_{(1)}\kappa_K) +
  \alpha_\theta,
\end{equation*}
with each $\kappa_A\in\Lambda^2_0E\subset S^2T^*M$ and $\alpha_\theta\in
S^2E$.  We have
\begin{gather*}
  \alpha_\theta = \tfrac14(1+I+J+K)F_\theta,\qquad \mu_A =
  \tfrac1n\Lambda_AF_\theta,\\
  A_{(1)}\kappa_A = -\tfrac14A(1+I+J+K)AF_\theta - \mu_A\omega_A.
\end{gather*}
We find that
\begin{equation*}
  \tfrac12(J+K)F_\theta = -\mu_I\omega_I - I_{(1)}\kappa_I + \alpha_\theta
  \in \Lambda^{1,1}_I. 
\end{equation*}

\begin{proposition}
  Suppose $W$ is obtained from an almost hyperHermitian manifold~$M$ by a
  twist of the $\mathbb R$-action generated by a symmetry~$X$ and using a
  curvature form $F_\theta$.  Then $W$ carries an almost hyperHermitian
  structure and the intrinsic almost quaternionic torsion is determined by
  the one-forms
  \begin{gather*}
    {\nu_4^I}^W = \nu_4^I + \tfrac2a\{ \mu_I(2n-1)IX^\flat - X\hook
    \alpha_\theta - IX\hook\kappa_I\},\\
    {\nu_3^I}^W = \nu_3^I +
    \tfrac4{a(2n+1)(n-1)}(nIX\hook\kappa_I-(n-1)X\hook\alpha_\theta).
  \end{gather*}
  and three-form components
  \begin{gather*}
    \begin{split}
      {\beta_I^{(3)}}^W \mkern-15mu &= \beta_I^{(3)} +
      \tfrac2{3a}(2X^\flat\wedge I_{(1)}\kappa_I + JX^\flat \wedge
      K_{(1)}\kappa_I - KX^\flat \wedge J_{(1)}\kappa_I)\eqbreak +
      \tfrac2{3a(n-1)}(2(X\hook\kappa_I)\wedge\omega_I -
      (KX\hook\kappa_I)\wedge\omega_J + (JX\hook\kappa_I)\wedge\omega_K),
    \end{split}\\
    \begin{split}
      {\beta_I^{(K)}}^W\mkern-15mu &= \beta_I^{(K)} - \tfrac
      2aX^\flat\wedge\alpha_\theta\smalleqbreak +\tfrac2{3a}(X^\flat\wedge
      I_{(1)}\kappa_I - JX^\flat\wedge K_{(1)}\kappa_I+ KX^\flat\wedge
      J_{(1)}\kappa_I)\smalleqbreak
      -\tfrac2{3a(2n+1)}((X\hook\kappa_I)\wedge\omega_I+
      (KX\hook\kappa_I)\wedge \omega_J +
      (JX\hook\kappa_I)\wedge\omega_K)\smalleqbreak
      -\tfrac2{a(2n+1)}((IX\hook\alpha_\theta)\wedge\omega_I +
      (JX\hook\alpha_\theta)\wedge\omega_J +
      (KX\hook\alpha_\theta)\wedge\omega_K).
    \end{split}
  \end{gather*}
\end{proposition}

\begin{proof}
  One first computes the following contraction formul\ae
  \begin{gather*}
    I\Lambda_I(X^\flat\wedge\alpha_\theta) = X\hook\alpha_\theta,\quad
    I\Lambda_I(\gamma\wedge\omega_I) = (2n-1)I\gamma,\quad
    I\Lambda_I(\gamma\wedge\omega_J) = J\gamma,\\
    I\Lambda_I(X^\flat\wedge I\kappa_A) = -IX\hook\kappa_A,\quad
    I\Lambda_I(X^\flat\wedge J\kappa_A) = JX\hook\kappa_A.
  \end{gather*}
  These lead directly to the claimed expressions for ${\nu_4^I}^W$ and
  ${\nu_3^I}^W$ and then give
  \begin{gather*}
    {\beta_{4I}}^W = \beta_{4I} + \tfrac2{a(2n-1)}(\mu_I(2n-1)X^\flat+
    IX\hook\alpha_\theta - X\hook\kappa_I) \wedge\omega_I,\\
    \begin{split}
      {\beta_{3I}^{(E)}}^W \mkern-15mu &= \beta_{3I}^{(E)} -
      \tfrac4{a(2n+1)(n-1)(2n-1)}(-nX\hook\kappa_I -
      (n-1)IX\hook\alpha_\theta)\wedge\omega_I\eqbreak -
      \tfrac2{a(2n+1)(n-1)}(-nKX\hook\kappa_I -
      (n-1)JX\hook\alpha_\theta)\wedge\omega_J\eqbreak -
      \tfrac2{a(2n+1)(n-1)}(-nJX\hook\kappa_I -
      (n-1)KX\hook\alpha_\theta)\wedge\omega_K.
    \end{split}
  \end{gather*}
  The remaining components of ${\beta_I}^W$ are then found
  using~\eqref{eq:beta-K3} with
  \begin{gather*}
    \mathcal L_I(\gamma\wedge\alpha_\theta) = \gamma\wedge\alpha_\theta,\quad
    \mathcal L_I(\gamma\wedge\omega_J) =
    2I\gamma\wedge\omega_K-\gamma\wedge\omega_J,\\
    \mathcal L_I(\gamma\wedge\omega_I) = \gamma\wedge\omega_I,\quad
    \mathcal L_I(\gamma\wedge J_{(1)}\kappa_A) = 2I\gamma\wedge
    K_{(1)}\kappa_A - \gamma\wedge J_{(1)}\kappa_A.
  \end{gather*}
\end{proof}

\begin{corollary}
  \label{cor:twisted-torsion}
  Twisting by $F_\theta\in \sym2E+\sym2H$ leaves $\xi_{33}$, $\xi_{3H}$ and
  $\xi_{K3}$ unchanged.  Furthermore,
  \begin{enumerate}
  \item\label{item:EH} if $3X^\theta = -2(n+2)\sum_A \mu_AAX^\flat$, then
    $\xi_{EH}$ is not affected;
  \item\label{item:S2E} if $F_\theta\in \sym2E$, then $\xi_{E3}$ is unaltered.
  \end{enumerate}
\end{corollary}

\begin{proof}
  The assumption on $F_\theta$ is equivalent to the vanishing of
  $\kappa_I$, $\kappa_J$ and $\kappa_K$.  We thus have
  ${\beta_I^{(3)}}^W=\beta_I^{(3)}$ and that
  ${\beta_I^{(K)}}^W-\beta_I^{(K)}$ is independent of~$I$, from which the
  invariance of the components in $(\Lambda^3_0E+K)\sym3E+\Lambda^3_0EH$
  follows.

  The change in $I\eta_I$ is $\tfrac{(n-1)}{2na}((2n+1)(\mu_IIX^\flat -
  \mu_JJX^\flat-\mu_KKX^\flat) - X\hook\alpha_\theta)$.  If
  $3X\hook\alpha_\theta+(2n+1)\sum_A\mu_AAX^\flat=0$, then there is no
  contribution to $\sum_A A\eta_A$ and hence no change in~$\xi_{EH}$.  On
  the other hand, if $F_\theta\in \sym2E$ then each $\mu_A=0$, and there is
  no contribution to $E\sym3H$.
\end{proof}

\begin{remark}
  Case~\itref{item:S2E} shows in particular that the QKT condition is
  preserved by twisting with $F_\theta\in\sym2E$.
\end{remark}

\section{Examples}
\label{sec:examples}

In this section we use the techniques developed in the previous sections to
compute the intrinsic torsion in a number of particular examples.  In the
first instance we consider examples which a almost hyperHermitian with each
Hermitian structure of the same Gray-Hervella type.  From
\cite{Cabrera-S:aHqg} we know that certain combinations can not occur.
Table~\ref{tab:aHall} gives an overview of which types may be obtained.
The $4n$-torus $T^{4n}=\mathbb H^n/\mathbb Z^{4n}$ is hyperK\"ahler and so
has intrinsic torsions~$0$.  The Hopf-like manifold $S^{4n-1}\times S^1 =
{\mathbb H^n\setminus\{0\}}/(q\mapsto 2q)$, is locally, but not globally
conformal, to the flat hyperK\"ahler metric, so each almost Hermitian
structure is of class $\Wc4$ and the almost quaternionic type is~$EH$.  The
other examples are described below.

\begin{table}[tp]
  \centering
  \begin{tabular}{ccccc}
    \toprule
    $I,J,K$ & $\{0\}$ & $\Wc1$ & $\Wc2$ & $\Wc{1+2}$ \\
    $\{0\}$ & $T^{4n}$ & impossible & impossible & impossible \\
    $\Wc3$ & $\HQ$ & $S^3{\times} T^9$ & $T^3{\times}
    M(k)^3$ & $T^3{\times} (\Gamma\backslash H)^3$ \\
    $\Wc4$ & $S^{4n-1}{\times} S^1$ & impossible & impossible & unknown \\
    $\Wc{3+4}$ & $S^3{\times} T^{4m+1}$ & $S^3\times T^9$ & $T^3{\times}
    (\Gamma\backslash H)^3$, $T^3{\times} M(k)^3$ &
    $T^3{\times}(\Gamma\backslash H)^3$ \\ 
    \bottomrule
  \end{tabular}
  \caption{Examples with common almost Hermitian structures}
  \label{tab:aHall}
\end{table}

\subsection{The manifold \texorpdfstring{$S^3 \times T^9$}{S3xT9}}

The sphere $S^3$ is isomorphic to the Lie group $\SP(1)$.  In the Lie
algebra $\sP(1)$ there is a basis $x$, $y$, $z$ such that $[x,y]=2z$,
$[z,x]=2y$ and $[y,z]=2x$.  From $x$, $y$, $z$ one can determine the left
invariant one-forms $a$, $b$, $c$ which constitute a basis for one-forms
and their exterior derivatives are given by $ d a = - 2 b \wedge c$, $d b =
- 2 c \wedge a$ and $d c = - 2 a \wedge b$.

In the product manifold $M = S^3\times T^9 $, we write $a_1$, $b_1$, $c_1$
for the one-forms corresponding to the factor $S^3$ and $a_i$, $b_i$,
$c_i$, $i=2,3,4$, for a basis of invariant one forms on~$T^9=(S^1)^9$.  On
$M$ we consider an almost hyperHermitian structure $I$, $J$, $K$ with
compatible metric $\inp\cdot\cdot = \sum_{i=1}^4 (a_i^2 + b_i^2 + c_i^2)$
and whose K\"ahler forms are given by
\begin{equation}
  \label{eq:omega-J}
  \omega_I  = a_2  a_1 + a_4  a_3 +b_2 
  b_1 + b_4  b_3 +c_2  c_1 + c_4  c_3,\detc(234),
\end{equation}
where $a_2a_1=a_2\wedge a_1$ and \enquote{etc.(234)} denotes the
corresponding equations obtained by simultaneously cyclically permuting
$(I,J,K)$ and $(2,3,4)$.  Their respective exterior derivatives and the
three-forms $\beta_A$ are given by
\begin{equation*}
  d \omega_I = 2 \sumcic_{abc} a_2  b_1  c_1, \quad
  \beta_I = - 2 \sumcic_{abc} \left( a_1  b_3  c_3 + a_1 
    b_4  c_4 \right), \detc(234).
\end{equation*}
Since $\Lambda_B d\omega_A = 0$, $A,B=I,J,K$, we have $\lambda_A = 0$ and
$\beta_A^{(E)} = 0$, for $A=I,J,K$.  Then, using
equation~\eqref{eq:beta-K3}, we obtain
\begin{gather*}
  \beta_I^{(K)} = - \tfrac23 \sumcic_{abc} \left(a_1 b_1 c_1 + a_1 b_2 c_2
    + a_1 b_3 c_3 + a_1 b_4 c_4 \right), \detc \\
  \beta_I^{(3)} = - \tfrac23 \sumcic_{abc} \left( -a_1 b_1 c_1 - a_1 b_2
    c_2 + 2 a_1 b_3 c_3 + 2 a_1 b_4 c_4 \right), \detc(234).
\end{gather*}
Using Proposition~\ref{prop:torsion-beta} we get $\psi^{(3)}_A = 0$,
$\psi^{(K)}_A = 0$, $ \theta^\xi = 0 = \theta^\xi_A$, for $A=I,J,K$,
$\psi^{(3)} \neq 0$ and $\psi^{(K)} \neq 0$.  Therefore, we have
\begin{equation*}
  \xi \in \Lambda^3_0 E \sym3 H + K H.
\end{equation*}
Note also that $\xi_{33} \neq 0$ and $\xi_{KH} \neq 0$.

Furthermore, if we consider the connection $\Nt = \LC + \tfrac12T$, where
$T$ is given by
\begin{equation*}
  \inp Y{T_X Z} = \tfrac16 \sum_{A=I,J,K}
  \bigl( \beta_A^{(K)} - \beta_A^{(3)} \bigr) (X,Y,Z),
\end{equation*}
we will obtain that this connection is metric and $\Nt I = \Nt J = \Nt K =
0$.  Thus, we have got an example of a quaternion-Hermitian manifold, which
is not QKT, admitting an $\SP(n)\SP(1)$-connection with skew-symmetric
$(0,3)$-torsion.

We compute the Nijenhuis $(0,3)$-tensors for $I$, $J$ and $K$ via
Proposition~\ref{prop:naddd}:
\begin{equation*}
  N_I = 2 a_1  b_1  c_1 - 2 \sumcic_{abc} a_1  b_2 
  c_2, \detc(234).
\end{equation*}
Since $N_I$, $N_J$, $N_K$ are skew-symmetric, then the almost Hermitian
structures are of a type that lies in~$\Wc{1+3+4}$.  However, $ A d^*
\omega_A = - \Lambda_A d\omega_A =0$ implies that such structures are of
type $\Wc{1+3}$.  The facts $4 d \omega_A \neq 3 AN_A$ and $N_A \neq 0$
respectively imply that the structures are not of the types $\Wc1$
and~$\Wc3$.

Finally, if we make a conformal change of metric $\inp\cdot\cdot^o = e^{2
\sigma} \inp\cdot\cdot$ as in Proposition~\ref{prop:conformal}, we obtain a
new quaternionic structure with
\begin{equation*}
  \xi^o \in \Lambda^3_0 E \sym3 H + K H + EH.
\end{equation*}
The new almost Hermitian structures are of type $\Wc{1+3+4}$.

\subsection{The manifold \texorpdfstring{$S^3 \times
T^{4m+1}$}{S3xT(4m+1)}}

In the product manifold $M = S^3 \times T^{4m+1}$, $m \geqslant 1$, we
write $n=m+1$, $a_2$, $a_3$, $a_4$ to denote the one-forms corresponding to
the factor $S^3$ and $a_1$, $a_i$, $i=5, \dots ,4n$, for a basis of
invariant one-forms on~$T^{4m+1}$.  On $M$ we consider an almost
hyperHermitian structure $I$, $J$, $K$ with compatible metric
$\inp\cdot\cdot = \sum_{i=1}^{4n} a_i \otimes a_i$ and whose K\"ahler forms
are given by the expressions
\begin{equation}
  \label{eq:omegaImany}
  \omega_I = \sum_{i=0}^{n-1} \left( a_{4i+2}  a_{4i+1} +
    a_{4i+4}  a_{4i+3}\right), \detc(234).
\end{equation}
Their respective exterior derivatives and the three-forms~$\beta_A$ are
given by
\begin{equation*}
  d \omega_I = - 2 a_1 a_3 a_4, \detc(234), \quad \beta_I = - 4 a_2 a_3
  a_4, \detc
\end{equation*}
Hence, by Proposition \ref{prop:HKT}, $M$~is an HKT-manifold.  Although
this example is already known, we wish to give a few more details.  We have
\begin{gather*}
  I d^* \omega_I = -2 a_1, \quad I \lambda_I = - \tfrac1n\, a_1, \quad
  \nu^I_3 = - \tfrac8{2n+1}\, a_1, \\ \nu^I_4 = - 4 a_1, \quad \theta^\xi =
  \theta_I^\xi = - \tfrac1{4(2n+1)}\, a_1, \detc,
\end{gather*}
so
\begin{equation*}
  \xi \in (K+E)H
\end{equation*}
with $\xi_{KH} \neq 0$ and $\xi_{EH} \neq 0$.  Also, the three almost
Hermitian structures are of type $\Wc{3+4} \setminus ( \Wc3 \cup \Wc4 )$.

Making a local conformal change of metric $\inp\cdot\cdot^o = e^{2 \sigma}
\inp\cdot\cdot$, with $\sigma$ satisfying $d \sigma = 4 \theta^\xi = -
\tfrac1{(2n+1)}\, a_1$, we will obtain a new almost hyperHermitian
structure $\{ I,J,K ; \inp\cdot\cdot^o \}$ with
\begin{equation*}
  \xi^o \in KH.
\end{equation*}
However, by Proposition \ref{prop:conformal}, we will have
\begin{gather*}
  I d^* \omega^o_I = - \tfrac4{(2n+1)}\, a_1, \quad I \lambda^o_I = -
  \tfrac2{2n+1} \, a_1, \quad {\nu^I_3}^o = - \tfrac4{2n+1}\, a_1, \\
  {\nu^I_4}^o = - \tfrac{8n}{2n+1} a_1, \quad \theta^{\xi^o} =
  \theta_I^{\xi^o} = 0, \detc,
\end{gather*}
so the identities given in Lemma~\ref{lem:hktonef} are not satisfied and
structure is not HKT.  It has $A \lambda^o_A = \tfrac1{2n} t^o \neq 0$, but
the three almost Hermitian structures are still of type $\Wc{3+4}\setminus(
\Wc3 \cup \Wc4)$.

Finally, for a local conformal change of metric by $\sigma$ satisfying $d
\sigma = -\tfrac1{2n-1} a_1$, we obtain an almost hyperHermitian structure
with three almost Hermitian structures of type~$\Wc3$.  However, such a
structure is not HKT.  The almost quaternion-Hermitian structure has
\begin{equation*}
  \xi^o\in(K+E)H 
\end{equation*}
but not in any submodule.  In fact, we will have
\begin{gather*}
  I d^* \omega^o_I = 0, \quad I \lambda^o_I = - \tfrac{2(n-1)}{2n-1} \,
  a_1, \quad {\nu^I_3}^o = - \tfrac{4(2n-3)}{4n^2-1}\, a_1, \\
  {\nu^I_4}^o = - \tfrac{8(n-1)}{2n-1} a_1, \quad \theta^{\xi^o} =
  \theta_I^{\xi^o} = \tfrac1{2(4n^2-1)}\, a_1, \detc
\end{gather*}

\subsection{The quaternionic Heisenberg group}

Now we consider the quaternionic Heisenberg group~$\HQ$ described
by Cordero et al.~\cite{Cordero-FL:Heisenberg}:
\begin{equation*}
  \HQ = \left\{
    \,
    \begin{pmatrix}
      1 & q_1 & q_3\\
      0 & 1 & q_2 \\
      0 & 0 & 1
    \end{pmatrix}
    : q_1,q_2,q_3\in\mathbb H
    \,
  \right\},
\end{equation*}
which is a connected nilpotent Lie group. A basis for the left-invariant
one-forms on~$\HQ$ is given by $a_i$, $b_i$, $c_i$, $i=1,\dots,4$, where
\begin{gather*}
  dq_1 = a_1 + ia_2 + ja_3 + ka_4,\quad
  dq_2 = b_1 + ib_2 + jb_3 + kb_4,\\
  dq_3 - q_1dq_2 = c_1+ ic_2 + jc_3 + kc_4.
\end{gather*}
With $\Gamma_Q$ be the subgroup of matrices of $\HQ$ with $q_i\in \mathbb
Z\{1,i,j,k\}$, we see that these forms descend to the compact manifold $M_Q
= \Gamma_Q \backslash \HQ$.  Consider the almost hyperHermitian structure
on $M_Q$ with metric $\inp\cdot\cdot = \sum_{i=1}^4 ( a_i^2 + b_i^2 +
c^2_i)$, and K\"ahler forms given by~\eqref{eq:omega-J}.

\begin{proposition}[Cordero et al.~\cite{Cordero-FL:Heisenberg}]
  The three almost Hermitian structures $I$, $J$, $K$ defined on $M_Q$ are
  of type~$\Wc3$. \qed
\end{proposition}

\noindent
Structures of type $\Wc3$ are sometimes called \emph{balanced} Hermitian.

The three-forms $\beta_A$ are given by 
\begin{equation*}
  \begin{split}
    \tfrac12 \beta_I &= a_1 b_1 c_1 + a_1 b_2 c_2 + a_1 b_3 c_3 + a_1 b_4
    c_4\eqbreak 
    - a_2 b_1 c_2 + a_2 b_2 c_1 - a_2 b_3 c_4 + a_2 b_4 c_3, \detc(234).
  \end{split}
\end{equation*}
Using equation~\eqref{eq:beta-theta}, we obtain that $\nu^A_3 = \nu^A_4 =
0$, so $\xi_{EH} = \xi_{E3} =0$ by Proposition~\ref{prop:torsion-beta}.

On the other hand, using equation~\eqref{eq:beta-K3}, we obtain
$\beta_I^{(K)}= \beta_J^{(K)} = \beta_K^{(K)} \neq 0$ and $\beta_I^{(3)}+
\beta_J^{(3)} + \beta_K^{(3)} = 0$, with each $\beta_A^{(3)} \neq 0$.
Therefore, Proposition~\ref{prop:torsion-beta} gives $\xi_{33} = \xi_{K3}
=0$, $\xi_{3 H} \neq 0$ and $\xi_{KH} \neq 0$.  In summary,
\begin{equation*}
  \xi \in \Lambda_0^3 EH + KH.
\end{equation*}

\subsection{The manifold \texorpdfstring{$T^3 \times
(\Gamma\backslash H)^3$}{T3x(H/Gamma)3}}

Let $H$ be the real Heisenberg group of dimension~$3$:
\begin{equation*}
  H = \left\{\,
    \begin{pmatrix}
      1 & x & z \\
      0 & 1 & y \\
      0 & 0 & 1
    \end{pmatrix}
    : x,y,z\in\mathbb R\,
  \right\}.
\end{equation*}
A basis of left-invariant one-forms is given by $\{ dx , dy , dz - x dy
\}$.  Let $\Gamma$ be the discrete subgroup of matrices of $H$ whose
entries $x$, $y$, $z$ are integers.  The quotient space $\Gamma \backslash
H $ is called the Heisenberg compact nilmanifold.  The one-forms $\{ dx ,
dy , dz - x dy \}$ descend to one-forms $p,q,r$ on $\Gamma\backslash H$
with $dr=-p\wedge q$.  

\subsubsection{A first  structure}

Let $M$ be the manifold $T^3\times (\Gamma\backslash H)^3$.  On $M$ we
consider a basis of invariant one-forms $a_1$, $b_1$, $c_1$ for~$T^3$ and
one-forms $a_2$, $a_3$, $a_4$, $b_2$, $b_3$, $b_4$, $c_2$, $c_3$ and $c_4$
corresponding to the factors $\Gamma\backslash H$ such that $ d a_2 = -
a_3 a_4$, $d b_3 = - b_4 b_2$ and $d c_4 = - c_2 c_3$.

On $M$ we consider an almost hyperHermitian structure $I$, $J$, $K$ with
compatible metric $\inp\cdot\cdot = \sum_{i=1}^4 (a_i^2 + b_i^2 + c_i^2)$
and whose K\"ahler forms are given by the expressions~\eqref{eq:omega-J}.

Their respective exterior derivatives and three-forms $\beta_A$ are given
by 
\begin{equation*}
  d \omega_I = - a_1 a_3 a_4,\quad \beta_I = - b_2 b_3 b_4 - c_2 c_3
  c_4, \detc(abc;234).
\end{equation*}
Using equations \eqref{eq:beta-E4}, \eqref{eq:beta-theta}
and~\eqref{eq:ladd2}, we obtain
\begin{equation*}
  \nu^I_3 = - \tfrac27 (b_1 +c_1) , \quad \nu^I_4 = -
  b_1 -c_1, \quad I \lambda_I = \tfrac16
  (a_1-b_1-c_1), \detc(abc).
\end{equation*}
Therefore, by Proposition \ref{prop:torsion-beta}, we get
\begin{gather*}
  56\,\theta^\xi_I = -3a_1+b_1+c_1, \quad -168\, \theta^\xi = a_1+b_1+c_1,
  \detc(abc).
\end{gather*}
Thus, the $\xi_{E3}$ and $\xi_{EH}$ parts of the intrinsic torsion are not
zero.

On the other hand, by equations \eqref{eq:beta-E3} and~\eqref{eq:beta-E4},
we have
\begin{equation*}
  7\beta_I^{(E)} = \sum_{A=I,J,K} A (b_1+c_1)  \omega_A, \detc(abc).
\end{equation*}
Hence, by equation~\eqref{eq:beta-K3} on $\beta_A^{(3+K)} = \beta_A -
\beta_A^{(E)}$, we have
\begin{gather*}
  \beta_I^{(3)} = 0, \quad \beta^{(K)}_I = - b_2 b_3 b_4 - c_2 c_3 c_4 -
  \tfrac17 \sum_{A=I,J,K} A (b_1+c_1) \omega_A, \detc(abc).
\end{gather*}
Proposition~\ref{prop:torsion-beta} then gives $\xi_{33} = \xi_{3H}= 0$,
$\xi_{K3} \neq 0$ and $\xi_{KH} \neq 0$ and we conclude
\begin{equation*}
  \xi \in K \sym3 H + E \sym3 H + KH + EH.
\end{equation*}

Let us now analyse the almost Hermitian structures.  We compute the
Nijenhuis tensors using Proposition~\ref{prop:naddd}:
\begin{equation*}
  \begin{split}
    N_I & =  b_3 \otimes b_1 b_3- b_3 \otimes b_2 b_4
    - b_4 \otimes  b_1  b_4 - b_4 \otimes  b_2  b_3\eqbreak
 - c_3 \otimes c_1 c_3 + c_3 \otimes c_2 c_4 + c_4 \otimes c_1 c_4 +
    c_4 \otimes c_2 c_3, \detc(abc;234).
  \end{split}
\end{equation*}
We see that the alternation of the Nijenhuis tensors $N_A$ are zero.
Therefore, the almost Hermitian structures are of type $\Wc{2+3+4}$.

Moreover, $N_A \neq 0$ so the structure is not of type~$\Wc{3+4}$.  Since
$I d^* \omega_I = - a_1$, etc.$(abc)$, the structures are not of
type~$\Wc{2+3}$.  It can be also checked that $5d \omega_A \neq - A d^*
\omega_A \omega_A$, so the structures also not of type $\Wc{2+4}$.

Finally, making a conformal change of metric $\inp\cdot\cdot^o = e^{2
\sigma} \inp\cdot\cdot$ with $\sigma$ is a local function satisfying
$d\sigma = 4 \theta^\xi = - \tfrac1{42} (a_1+b_1+c_1)$, we obtain locally a
new almost quaternion-Hermitian structure $\{ I,J,K ; \inp\cdot\cdot^o \}$
with
\begin{equation*}
  \xi^o \in K \sym3 H + E \sym3 H + KH.
\end{equation*}
The three new almost
Hermitian structures are still in $\Wc{2+3+4}$.

\subsubsection{A second  structure}

This time, on $M = T^3\times(\Gamma\backslash H)^3$, consider $a_1$, $b_1$,
$c_1$ a basis of invariant one-forms on~$T^3$, as before, and let $a_2$,
$a_3$, $a_4$, $b_2$, $b_3$, $b_4$, $c_2$, $c_3$, $c_4$ denote linearly
independent one-forms on the factors $\Gamma\backslash H$ now with $d a_2 =
- b_2 c_2$, $d b_3 = - c_3 a_3$ and $d c_4 = - a_4 b_4$.

Take the almost hyperHermitian structure $I$, $J$, $K$ with compatible the
metric $\inp\cdot\cdot = \sum_{i=1}^4 (a_i^2 + b_i^2 + c_i^2)$ and with
K\"ahler forms are respectively given by~\eqref{eq:omega-J}.

Their respective exterior derivatives and three-forms $\beta_A$ are given
by
\begin{gather*}
  d \omega_I = - a_1 b_2 c_2 + a_3 b_4 c_3 - a_4
  b_4  c_3, \detc(abc;234),\\
  \beta_I = - a_1 b_1 c_4 - a_1 b_3 c_1 - a_2 b_2 c_4 - a_2 b_3 c_2 - a_2
  b_3 c_3 - a_2 b_4 c_4, \detc(abc;234).
\end{gather*}
Note that $\Lambda_B d\omega_A = 0$ and $ \Lambda_B \beta_A = 0$, for
$A,B=I,J,K$, so $A \lambda_A=0$ and $\beta_A^{(E)}=0$.  This implies that
$\xi_{E3} = \xi_{EH} =0$.  Furthermore, using equation~\eqref{eq:beta-K3}
one computes the components $\beta_A^{(3)}$ and $\beta_A^{(K)}$ and obtains
that $\sum_{A=I,J,K} \beta_A^{(K)} \neq 0$ and $ 3 \beta_A \neq
\sum_{B=I,J,K} \beta_B^{(K)}$.  By Proposition~\ref{prop:torsion-beta}, we
find that $\xi_{K3} \neq 0$ and $\xi_{KH} \neq 0$.  Also we compute
$\sum_{A=I,J,K} \beta_A^{(3)} \neq 0$ and $6 \beta_A^{(3)} \neq 3
\sum_{B=I,J,K} \beta_B^{(3)} + \mathcal L_A \left( \sum_{B=I,J,K}
  \beta_B^{(3)} \right)$, so $\xi_{33} \neq 0$ and $\xi_{3H} \neq 0$.  In
summary, we get
\begin{equation*}
  \xi \in \Lambda_0^3 E \sym3 H + K \sym3 H + KH + EH.
\end{equation*}

Now we turn to analysis of the almost Hermitian structures.  Since $A d^*
\omega_A = - \Lambda_A d\omega_A = 0$, the almost Hermitian structures are
of type $\Wc{1+2+3}$.  Moreover, using the expressions for~$N_A$ in
Proposition~\ref{prop:naddd}, we have
\begin{equation*}
  \begin{split}
    N_I & = - a_1 \otimes b_1 c_2 -a_1 \otimes b_2 c_1 + a_2 \otimes b_2
    c_2 - a_2 \otimes b_1 c_1 \eqbreak + b_3 \otimes c_3 a_3 - b_3 \otimes
    c_4 a_4 - b_4 \otimes c_3 a_4 - b_4 \otimes c_4 a_3 \eqbreak - c_3
    \otimes a_3 b_3 - c_3 \otimes a_4 b_3 - c_4
    \otimes a_3  b_3 + c_4 \otimes a_4  b_4, \detc(abc;234).
  \end{split}
\end{equation*}
We find that $N_A$ is not skew-symmetric, the alternation ${\mathcal N}_A$
of $N_A$ is non-zero and that $4 d \omega_A \neq - A {\mathcal N}_A$.
Therefore, the almost Hermitian structures are not of type $\Wc{i+j}$, for
$i,j=1,2,3$.

Making a global conformal change we obtain almost Hermitian structures of
the general type $\Wc{1+2+3+4}$ whilst preserving the almost quaternion
Hermitian type.

\subsection{The manifold \texorpdfstring{$T^3\times M(k)^3$}{T3xM(k)3}}

Let us consider the manifolds $M(k)$ described in \cite{Auslander-GH:flows}
as $S_1/D$ and studied geometrically in
\cite{Cordero-FG:minimales,Cordero-FdLS:symplectic-4,Fernandez-G:symplectic-solv}.
For a fixed $k\in \mathbb R\setminus\{0\}$, let $\Gk$ be the
three-dimensional connected and solvable (non-nilpotent) Lie group
consisting of matrices:
\begin{equation*}
  \Gk =
  \left\{\,
    \begin{pmatrix}
      e^{kz} & 0       & 0 & x \\
      0      & e^{-kz} & 0 & y \\
      0      & 0       & 1 & z \\
      0 & 0 & 0 & 1
    \end{pmatrix}
    : x,y,z\in\mathbb R
    \,
  \right\}.
\end{equation*}
A basis of right invariant one-forms on $\Gk$ is $\{ dx - kx\,dz , dy +
ky\,dz, dz \}$.

The Lie group $\Gk$ possesses a discrete subgroup $\Gamma(k)$ such that the
manifold $M(k) = \Gk / \Gamma(k)$ is compact.  One example of $\Gamma(k)$
is generated by choosing $r>0$ so that $e^{kr}+e^{-kr}\in \mathbb Z$ and
then taking the subgroup of~$\Gk$ generated by $(x,y\in \mathbb Z,z=0)$ and
$(x=0=y,z=r)$.  The given basis of one-forms on $\Gk$ descends to one-forms
$a,b,c$ on~$M(k)$ satisfying $ da = - k ac$ and $db = k bc$.

\subsubsection{A first  structure}

Let $M$ be the manifold $M = T^3 \times M(k)^3$.  Consider a basis of
invariant one-forms $a_1$, $b_1$, $c_1$ on~$T^3$ and one-forms $a_2$,
$a_3$, $a_4$, $b_2$, $b_3$, $b_4$, $c_2$, $c_3$, $c_4$ corresponding on the
factors $M(k)$ with
\begin{equation*}
  d a_3 = - k a_3  a_2, \quad d a_4 = k a_4  a_2, \detc(abc;234).
\end{equation*}
Now we consider on $M$ an almost hyperHermitian structure $I$, $J$, $K$
with compatible metric $\inp\cdot\cdot = \sum_{i=1}^4 (a_i^2 + b_i^2 +
c_i^2)$ and K\"ahler forms given by~\eqref{eq:omega-J}.

The respective exterior derivatives and three-forms~$\beta_A$ are given by
\begin{equation*}
  d \omega_I = k b_1  b_2  b_3 - k c_1  c_2  c_4 ,
  \quad \beta_I = - k b_1  b_2  b_4 - k c_1  c_2 c_3,\detc(abc;234).
\end{equation*}
Note that $\sum_A \beta_A = 0$.  Therefore, $\sum_A \beta_A^{(3)} = 0$ and
$\sum_A \beta_A^{(K)} = 0$.  This implies that $\xi_{33}=0$ and $\xi_{K3}
=0$.  To compute the $E$-parts of $\beta_A$ and $\xi$, we first find
\begin{gather*}
  I d^* \omega_I = k b_3 - k c_4, \quad
  I\Lambda_K d\omega_J = - k a_2 + k c_4, \\
  I\Lambda_J d\omega_K = - k a_2 + k b_3, \detc(abc;234).
\end{gather*}
Using equation~\eqref{eq:ladd2}, we now obtain $I \lambda_I = \tfrac{k}3
(-b_3+c_4)$, etc.(abc;234).  On the other hand, using $\nu^A_4 =
A\Lambda_A\beta_A$ and equation~\eqref{eq:beta-theta}, we get
\begin{equation*}
  \nu^I_4= k (-b_3+c_4), \quad
  \nu^I_3= \tfrac{2k}7 (-b_3+c_4), \detc(abc;234).
\end{equation*}
Finally, from the obtained expressions for $A\lambda_A$, $\nu^A_3$ and
$\nu^I_4$, using Proposition~\ref{prop:torsion-beta}, we find
\begin{gather*}
  \theta^\xi = 0, \quad \theta_I^\xi= \tfrac{5k}{28} (-b_3+c_4),
  \detc(abc;234).
\end{gather*}
Hence $\xi_{EH}=0$ and $\xi_{E3} \neq 0$.  Furthermore, since
\begin{equation*}
  \beta_I^{(E)} = \tfrac{k}7 \sum_{A=I,J,K} A (b_3 - c_4) 
  \wedge \omega_A, \detc(abc;234),
\end{equation*}
we have $\mathcal L_A ( \beta_B - \beta_{B}^{(E)} ) = \beta_B -
\beta_{B}^{(E)}$, for $A,B=I,J,K$.  Therefore, $ \beta_A^{(K)} = \beta_A -
\beta_A^{(E)} \neq 0$ and $\beta_I^{(3)} = \beta_J^{(3)} = \beta_K^{(3)}
=0$.  These last claims imply $\xi_{K3} \neq 0$ and $\xi_{3H} =0$.  In
conclusion,
\begin{equation*}
  \xi \in K \sym3 H + E \sym3 H.
\end{equation*}

Now we analyse the almost Hermitian structures.  Using the expression for
$N_A$ in Proposition \ref{prop:naddd}, one can obtain Nijenhuis
$(0,3)$-tensor for $I$, $J$ and $K$.  These tensors $N_A$ are not zero and
but their alternations $\mathcal N_A$ are zero.  Therefore, the almost
Hermitian structures are of type $\Wc{2+3+4}$.

Moreover, $N_A \neq 0$, this implies that the structures are not of type
$\Wc{3+4}$.  Additionally, we have seen above that the Lee one-forms $A d^*
\omega_A$ are non-zero, so the structures are not of type~$\Wc{2+3}$.
Finally, one can easily check that $d \omega_A \neq - \tfrac15 A d^*
\omega_A \omega_A$, so the almost Hermitian structures are not of
type~$\Wc{2+4}$.

\subsubsection{A second  structure}

We again consider $M = T^3 \times M(k)^3$.  We take $a_1$, $b_1$, $c_1$ to
be a basis of invariant one-forms on~$T^3$ and now $a_2$, $a_3$, $a_4$,
$b_2$, $b_3$, $b_4$, $c_2$, $c_3$, $c_4$ is basis of one-forms on $M(k)^3$
with
\begin{equation*}
  d b_2 = - k b_2  a_2, \quad d c_2 = k c_2  a_2, \detc(abc;234).
\end{equation*}
On $M$ we consider an almost hyperHermitian structure $I$, $J$, $K$ with
compatible metric $\inp\cdot\cdot = \sum_{i=1}^4 (a_i^2 + b_i^2 + c_i^2)$
and K\"ahler forms given by~\eqref{eq:omega-J}.  The respective exterior
derivatives and the three-forms $\beta_A$ are given by
\begin{gather*}
  \tfrac1k \, d \omega_I = a_3 a_4 b_3 - a_3 a_4 c_4 - b_1 b_2 a_2 + b_3
  b_4 c_4 + c_1  c_2  a_2 -   c_3  c_4  b_3,\\
  \begin{split}
    \tfrac1{k} \, \beta_I & = - a_1 a_3 b_1 + a_1 a_4 c_1 - a_2 a_3 b_2 +
    a_2 a_4 c_2 - b_1 b_4 c_1 - b_2 b_3 a_3 \eqbreak - b_2 b_4 a_4 - b_2
    b_4 c_2 + c_1 c_3 b_1 + c_2 c_3 a_3 + c_2 c_3 b_2 + c_2 c_4 a_4,
  \end{split}
\end{gather*}
etc.(abc;234).  Since $\Lambda_B d \omega_A = 0$, for $A,B=I,J,K$, we have
$ A \lambda_A = \nu^A_3 = \nu^A_4 = \theta^\xi = \theta_A^\xi = 0$, so
$\xi_{E3} =0$ and $ \xi_{EH}=0$.  Now, using~\eqref{eq:beta-K3}, one can
compute the components $\beta_A^{(3)}$ and $\beta_I^{(K)}$.  One checks
that $\sum_{A=I,J,K} \beta_A^{(3)}=0$, $\beta_A^{(3)}\neq 0$, the
$\beta_A^{(K)}$, $A=I,J,K$, are distinct and $\sum_{A=I,J,K} \beta_A^{(K)}
\neq 0$.  Therefore, $\xi_{33} = 0$, $\xi_{EH} \neq 0$, $\xi_{K3}\neq 0$
and $\xi_{KH} \neq 0$.  In summary,
\begin{equation*}
  \xi \in K \sym3 H + \Lambda_0^3 EH + KH.
\end{equation*}
 
For the almost Hermitian structures, one computes $N_A$ via
Proposition~\ref{prop:naddd} and find that they are non-zero but their
alternations ${\mathcal N}_A$ vanish, so the structures lie in
$\Wc{2+3+4}$.  Since $\Lambda_B d \omega_A = 0$, for $A,B=I,J,K$, we have
that the respective $\Wc4$-parts are zero.  Thus, the almost Hermitian
structures are of type $\Wc{2+3}$.  Because $N_A\neq 0$ and $d \omega_A
\neq 0 $, the structures are not of the simple types $\Wc2$ or~$\Wc3$.
 
\subsection{Twisting tori}
\label{sec:tori}

Up to this point we have obtained all the claimed examples of
Table~\ref{tab:aHall} and determined the corresponding types of the almost
quaternion-Hermitian structures.  Let us now use the twist construction
of~\S\ref{sec:twist} to give some other examples of types of~$\xi$.

Let us first demonstrate that the condition in
Corollary~\ref{cor:twisted-torsion}\itref{item:EH} can be satisfied.  Let
$M=T^{4n}$, $n\geqslant2$, with the standard flat hyperK\"ahler structure.
Take $F_\theta = \omega_I+\alpha_\theta$ and $X$~an arbitrary isometry.
Then $\mu_I=1$, $\mu_J=0=\mu_K$ and $X^\theta = -IX^\flat +
X\hook\alpha_\theta$.  The condition of
Corollary~\ref{cor:twisted-torsion}\itref{item:EH} is now equivalent to
$X\hook\alpha_\theta = -\tfrac{(2n+1)}3IX^\flat$.  This is satisfied if we
take $\alpha_\theta = -\tfrac{(2n+1)}{3\norm X^2}(X^\flat\wedge IX^\flat -
JX^\flat\wedge KX^\flat)$.  Locally, we may now twist to obtain a structure
with
\begin{equation*}
  \xi^W \in KH + ES^3H.
\end{equation*}
Since $X^\flat\wedge\alpha_\theta\ne0$, $\alpha$ is supported on~$\mathbb
HX$ and $n\geqslant2$, we have ${\xi_{KH}}^W\ne0$.  On the other hand the
values of the $\mu_A$ ensure that ${\xi_{E3}}^W\ne0$.

As second example, consider $T^{4n}$, $n\geqslant3$ with invariant basis
$a_1,\dots,a_{4n}$ and two-forms~\eqref{eq:omegaImany}.  Take $F_\theta =
a_2a_1+a_4a_3-a_6a_5-a_8a_7$.  We have $F_\theta\in \Lambda^{1,1}_I$
orthogonal to $\omega_I$ and of type $\{2,0\}$ for $J$ and $K$, so
$F_\theta = I_{(1)}\kappa_I$.  Twisting via any $X$ such that $X\hook
F_\theta=0$ and each $AX\hook F_\theta=0$ yields an almost
quaternion-Hermitian structure with
\begin{equation*}
  \xi^W \in (\Lambda^3_0E+K)(S^3H+H)
\end{equation*}
but not in any proper submodule.  Making a conformal change we may also
obtain structures with intrinsic torsion in
$(\Lambda^3_0E+K)S^3H+(\Lambda^3_0E+K+E)H$.

\subsection{Twisting Salamon's example}
\label{sec:Salamon}

Recall that Salamon~\cite{Salamon:parallel} gave an example of a
non-quaternionic K\"ahler $8$-manifold that has $d\Omega=0$.  The manifold
is a compact nilmanifold $\Gamma\backslash G$, where $G$ has a basis
of left-invariant one-forms $a_1,\dots,a_4,b_1,\dots,b_4$ satisfying
\begin{gather*}
  da_2 = -\sqrt3 a_1a_4 - 3a_4b_1,\quad
  db_2 = a_1a_4 + \sqrt3 a_4b_1,\\
  db_4 = - a_1a_2 - \sqrt3a_2b_1 - \sqrt3a_1b_2 + 3b_1b_2,
\end{gather*}
with the other basis elements closed.  The Lie algebra $\lie g$ is a direct
sum $\mathbb R^3 + \lie h$, with $\lie h$ two-step nilpotent.  Salamon's
structure is then given by \eqref{eq:omega-J} (with $c_i=0$).  In terms of
intrinsic torsion this has
\begin{equation*}
  \xi \in K\sym3H.
\end{equation*}

The two-form $F_\theta = a_1a_3+a_2a_4$ is a closed element of $S^2E$ and
defines an integral cohomology class on~$M$ for an appropriate choice
of~$\Gamma$.  We may thus twist using, for example, the central vector
field~$X$ dual to~$b_3$ to obtain an almost quaternion Hermitian
$8$-manifold~$W$ with
\begin{equation*}
  \xi^W \in K(\sym3H+H).
\end{equation*}

Conformally scaling these two examples we may obtain structures with
\begin{equation*}
  \xi^o \in K\sym3H + EH\quad\text{and}\quad {\xi^W}^o \in K\sym3H +
  (K+E)H. 
\end{equation*}


\begin{thebibliography}{10}

\bibitem{Auslander-GH:flows}
L.~Auslander, L.~Green, and F.~Hahn, \emph{Flows on homogeneous spaces}, With
  the assistance of L. Markus and W. Massey, and an appendix by L. Greenberg.
  Annals of Mathematics Studies, No. 53, Princeton University Press, Princeton,
  N.J., 1963.

\bibitem{Berger:hol}
M.~Berger, \emph{Sur les groupes d'holonomie des vari{\'e}t{\'e}s {\`a}
  connexion affine et des vari{\'e}t{\'e}s {r}iemanniennes}, Bull. Soc. Math.
  France \textbf{83} (1955), 279--330.

\bibitem{Broecker-tom-Dieck:Lie}
{Th}. Br{\"o}cker and T.~tom Dieck, \emph{Representations of compact {L}ie
  groups}, Graduate Texts in Mathematics, vol.~98, Springer-Verlag, New York,
  1985.

\bibitem{Cleyton-S:intrinsic}
R.~Cleyton and A.~F. Swann, \emph{{E}instein metrics via intrinsic or parallel
  torsion}, Math. Z. \textbf{247} (2004), no.~3, 513--528.

\bibitem{Cordero-FL:Heisenberg}
L.~A. Cordero, M.~Fern{\'a}ndez, and M.~de~Le{\'o}n, \emph{On the quaternionic
  {H}eisenberg group}, Boll. Un. Mat. Ital. A (7) \textbf{1} (1987), no.~1,
  31--37.

\bibitem{Cordero-FG:minimales}
L.~A. Cordero, M.~Fern{\'a}ndez, and A.~Gray, \emph{Modelos minimales en
  geometri{\'\i}a diferencial (minimal models in differential geometry)},
  Proceedings of the Workshop on Recent Topics in Differential Geometry (Puerto
  de la Cruz, 1990) (La Laguna), Informes, vol.~32, Univ. La Laguna, 1991,
  pp.~31--41. \MR{MR1127451 (92h:53038)}

\bibitem{Cordero-FdLS:symplectic-4}
L.~A. Cordero, M.~Fern{\'a}ndez, M.~de~Le{\'o}n, and M.~Saralegui,
\emph{Compact symplectic four solvmanifolds without 
  polarizations}, Ann. Fac. Sci. Toulouse Math. (5) \textbf{10} (1989), no.~2,
  193--198.

\bibitem{Fernandez-Moreiras:symmetry}
M.~Fern{\'a}ndez and B.~R. Moreiras, \emph{Symmetry properties of the covariant
  derivative of the fundamental {$4$}-form of a quaternionic manifold}, Riv.
  Mat. Univ. Parma \textbf{12} (1986), 249--256.

\bibitem{Fernandez-G:symplectic-solv}
M.~Fern{\'a}ndez and A.~Gray, \emph{Compact symplectic solvmanifolds
  not admitting complex structures}, Geom. Dedicata \textbf{34} (1990), no.~3,
  295--299.

\bibitem{Gates-HR:twisted}
S.~J. Gates, Jr., C.~M. Hull, and M.~Ro{\v c}ek, \emph{Twisted multiplets and
  new supersymmetric non-linear {$\sigma$}-models}, Nucl. Phys.~B \textbf{248}
  (1984), 157--186.

\bibitem{Grantcharov-P:HKT}
G.~Grantcharov and Y.~S. Poon, \emph{Geometry of hyper-{K\"a}hler connections
  with torsion}, Comm. Math. Phys. \textbf{213} (2000), no.~1, 19--37.

\bibitem{Gray:minimal}
A.~Gray, \emph{Minimal varieties and almost {H}ermitian submanifolds}, Michigan
  Math. J. \textbf{12} (1965), 273--287.

\bibitem{Gray-H:16}
A.~Gray and L.~M. Hervella, \emph{The sixteen classes of almost {H}ermitian
  manifolds and their linear invariants}, Ann. Mat. Pura Appl. (4) \textbf{123}
  (1980), 35--58.

\bibitem{Hitchin:Riemann-surface}
N.~J. Hitchin, \emph{The self-duality equations on a {R}iemann surface}, Proc.
  London Math. Soc. \textbf{55} (1987), 59--126.

\bibitem{Howe-OP:QKT}
P.~S. Howe, A.~Opfermann, and G.~Papadopoulos, \emph{Twistor spaces for {QKT}
  manifolds}, Comm. Math. Phys. \textbf{197} (1998), no.~3, 713--727.

\bibitem{Howe-P:further}
P.~S. Howe and G.~Papadopoulos, \emph{Further remarks on the geometry of
  two-dimensional nonlinear {$\sigma$} models}, Classical Quantum Gravity
  \textbf{5} (1988), no.~12, 1647--1661.

\bibitem{Howe-P:twistor-kaehler}
\bysame, \emph{Twistor spaces for hyper-{K}\"ahler manifolds with torsion},
  Phys. Lett. B \textbf{379} (1996), no.~1-4, 80--86.

\bibitem{Ivanov:QKT}
S.~Ivanov, \emph{Geometry of quaternionic {K\"a}hler connections with torsion},
  J. Geom. Phys. \textbf{41} (2002), no.~3, 235--257.

\bibitem{Joyce:hypercomplex}
D.~Joyce, \emph{Compact hypercomplex and quaternionic manifolds},
  J.~Differential Geom. \textbf{35} (1992), 743--761.

\bibitem{Cabrera-S:aHqg}
F.~Mart{\'\i}n~Cabrera and A.~F. Swann, \emph{Almost {H}ermitian structures and
  quaternionic geometries}, Differential Geom. Appl. \textbf{21} (2004), no.~2,
  199--214.

\bibitem{Cabrera:aqh}
F.~Mart{\'{\i}}n~Cabrera, \emph{Almost quaternion-{H}ermitian
  manifolds}, Ann. Global Anal. Geom. \textbf{25} (2004), no.~3, 277--301.

\bibitem{Obata:connection}
M.~Obata, \emph{Affine connections on manifolds with almost complex, quaternion
  or {H}ermitian structure}, Japan J. Math. \textbf{26} (1956), 43--77.

\bibitem{Salamon:Invent}
S.~M. Salamon, \emph{Quaternionic {K\"a}hler manifolds}, Invent. Math.
  \textbf{67} (1982), 143--171.

\bibitem{Salamon:quaternionic}
\bysame, \emph{Differential geometry of quaternionic manifolds}, Ann. Scient.
  {\'E}c. Norm. Sup. \textbf{19} (1986), 31--55.

\bibitem{Salamon:holonomy}
\bysame, \emph{Riemannian geometry and holonomy groups}, Pitman Research Notes
  in Mathematics, vol. 201, Longman, Harlow, 1989.

\bibitem{Salamon:parallel}
\bysame, \emph{Almost parallel structures}, Global differential geometry:
  the mathematical legacy of Alfred Gray (Bilbao, 2000), Contemp. Math., vol.
  288, Amer. Math. Soc., Providence, RI, 2001, pp.~162--181.

\bibitem{Swann:symplectiques}
A.~F. Swann, \emph{Aspects symplectiques de la g{\'e}om{\'e}trie
  quaternionique}, C.~R. Acad. Sci. Paris \textbf{308} (1989), 225--228.

\bibitem{Swann:T}
\bysame, \emph{{T} is for twist}, preprint PP-2006-23, Department of
  Mathematics and Computer Science, University of Southern Denmark, July 2006,
  Proceedings of the XV International Workshop on Geometry and Physics, Spanish
  Royal Mathematical Society, to appear.

\end{thebibliography}
\providecommand{\bysame}{\leavevmode\hbox to3em{\hrulefill}\thinspace}
\providecommand{\MR}{\relax\ifhmode\unskip\space\fi MR }
\providecommand{\MRhref}[2]{%
  \href{http://www.ams.org/mathscinet-getitem?mr=#1}{#2}
}
\providecommand{\href}[2]{#2}

\end{document}